\documentclass[12pt]{article}
\usepackage[cp1251]{inputenc}
\usepackage{amsmath}
\usepackage{amssymb}
\usepackage{amscd}
\usepackage[russian]{babel}
\begin{document}
\author{Белошапка В.К.}

\date{13.12.2021}

\title{{\bf   Модельные $CR$-поверхности: \\ взвешенный подход}}

\maketitle

\begin{abstract}
В работе дано систематическое построение теории "взвешенных"$ \,$ модельных поверхностей, основанной на понятии типа CR-многообразия по Блуму-Грэму-Степановой. Основой построения является конструкция Пуанкаре. Показано, как использование взвешенных модельных поверхностей расширяет возможности метода. Ставятся новые вопросы.
\end{abstract}

\footnote{
Механико-математический факультет Московского университета им.Ломоносова,
Воробьевы горы, 119992 Москва, Россия, vkb@strogino.ru.\\
Московский центр фундаментальной и прикладной математики МГУ}

\begin{center}
{\bf 1.Введение}
\end{center}

\vspace{3ex}

Среди подходов к локальному анализу $CR$-подмногообразий комплексного пространства следует отметить метод модельной поверхности.
Это эффективный аналитический подход с более чем столетней историей.  Основой этого метода является некая версия теоремы о неявном отображении в классе формальных степенных рядов. Автором этого приема является А.Пуанкаре, который применял его как в задачах небесной механики, так и в той области, которая сейчас называется $CR$-геометрией \cite{P}.

Недавно сфера применимости этого подхода была расширена до класса произвольных ростков конечного типа по Блуму-Грэму \cite{VB20}. Однако при этом использовалась стандартная версия данного метода.  Речь идет о том, что все переменные из комплексной касательной ростка имеют одинаковый вес.  В ряде частных случаев давно и успешно использовалась более гибкая техника свободного назначения весов переменным комплексной касательной \cite{VB96},\cite{AE01},\cite{KMZ14}. Поясним сказанное примером  из работы \cite{AL97}.

\vspace{3ex}

  Рассмотрим гиперповерхность в пространстве  ${\bf C}^{n+2}$  с координатами $(z_1,...,z_n,\,\zeta, \, w=u+i\,v)$, заданную уравнением
  \begin{equation}\label{AL}
  v = 2 \; {\rm Re} (z_1 \bar{\zeta} + \dots + z_n \bar{\zeta}^n).
  \end{equation}
  Данная гиперповерхность имеет тип по Блуму-Грэму $m=(2)$, причем модельная поверхность $Q=\{v=2 \; {\rm Re} (z_1 \bar{\zeta})\}$ голоморфно вырождена и размерность алгебры ее автоморфизмов бесконечна, в то время как алгебра автоморфизмов самой гиперповерхности конечномерна.  Это то, что мы получаем при стандартном подходе, когда веса всех координат из комплексной касательной  -- и $z$, и $\zeta$ -- одинаковы и равны 1. Но если расставить веса иначе, а именно: $([\zeta]=1, \; [z_j]=1+n-j,\; [w]=n+1)$, то поверхность становится взвешенно однородной (веса $n+1$).  К тому же эта гиперповерхность голоморфно однородна.  Используя $\,$ "взвешенную" $\,$ версию конструкции Пуанкаре,  можно получить оценку размерности алгебры автоморфизмов для ростка возмущенной гиперповерхности
 $$v = 2 {\rm Re} \,(z_1 \bar{\zeta} + \dots + z_n \bar{\zeta}^n) +o(n+1).$$
Ясно, что правильная точка зрения на этот пример -- это рассмотрение  этой гиперповерхности как взвешенно однородной и что за этим примером  стоит  "взвешенная" $\,$ теория модельных поверхностей конечного типа. При этом первый шаг в построении такой теории, а именно доказательство "взвешенного" $\,$  аналога теоремы Блума-Грэма, уже сделан М.Степановой в работе \cite{MS2}.

\vspace{3ex}

Целью данной работы  является систематическое построение такой "взвешенной" $\,$ теории модельных поверхностей,  которая опирается на взвешенный тип ростка по Блуму-Грэму-Степановой  так же, как теория, построенная в работе \cite{VB20}, --  на тип по Блуму-Грэму \cite{BG}.

\vspace{3ex}

{\bf 2.Веса и  взвешенный тип ростка CR-многообразия}

\vspace{3ex}
Конструкция Пуанкаре, которая лежит в основе метода, является по преимуществу аналитической. Поэтому из двух эквивалентных определений  типа ростка для нас главным является аналитическое определение с помощью вида локальных уравнений ростка многообразия.

Отличие взвешенной версии от стандартной заключается в возможности выбора {\it разных}  весов для координат комплексной касательной.

Пусть координаты объемлющего пространства ${\bf C}^N$ поделены на две группы $z \in \mathbf{C}^n$ и $w \in \mathbf{C}^K$.
Соответственно каждое из двух координатных подпространств разложено в прямую сумму
$$\mathbf{C}^n=\mathbf{C}^{n_1}+\dots+\mathbf{C}^{n_{p}}, \quad \mathbf{C}^K=\mathbf{C}^{k_1}+\dots+\mathbf{C}^{k_{q}}.$$
Пусть задано два набора натуральных чисел
$$ \mu_1 < \mu_2 < \dots < \mu_{p},  \quad   m_1 < m_2 < \dots < m_{q}.$$
Назначим веса координатам следующим образом
$$[z_{\alpha}]=\mu_{\alpha}, \; \mbox{где  } z_{\alpha} \in \mathbf{C}^{n_{\alpha}}, \quad [w_{\beta}]=m_{\beta}, \; \mbox{где  } w_{\beta} \in \mathbf{C}^{k_{\beta}}.$$
Подразумевается, что
$$\bar{z}_{\alpha}, \;  \bar{w}_{\beta},  \; {\rm Re} \,z_{\alpha},  \;{\rm Im} \,z_{\alpha},  \; {\rm Re} \,w_{\beta}=u_{\beta}, \; {\rm Im} \,w_{\beta}= v_{\beta}$$
имеют соответствующие веса.  Причем дифференцирование по координате веса $\nu$ получает вес $-\nu$.

Пусть локальные уравнения ростка $CR$-многообразия $M_0$ в начале координат имеют вид
\begin{equation}\label{1}
v_{\beta}=\Phi_{\beta} (z,\bar{z},u) +o(m_{\beta}),   \quad   \beta=1, \dots, q,
\end{equation}
где $\Phi_{\beta}$ -- квазиоднородный вектор-полином веса $m_{\beta}$, а под $o(\nu)$ понимается вектор-функция, чье тейлоровское разложение в начале координат содержит члены веса строго больше, чем $\nu$.

В \cite{MS} показано, что несложными обратимыми полиномиальными заменами координат в окрестности начала локальные уравнения вида (\ref{1}) всегда можно подвергнуть дополнительной редукции и привести к такому же виду, где взвешенно однородные формы $\Phi_j$ удовлетворяют дополнительным условиям. Эти условия связаны с существованием обратимых голоморфных преобразований, сохраняющих вид (\ref{1}), но изменяющих младшие координатные формы уравнения.  Каждый голоморфный полином в пространстве $\mathbf{C}^N$ -- это линейная комбинация  мономов вида (используем мульти-обозначения) $ z_1^{\gamma_1} \dots z_p^{\gamma_p} \,
w_{1}^{\delta_{1}} \dots  w_{q}^{\delta_{q}} $.  Каждый такой моном -- это возможность менять некоторые младшие координатные формы в уравнении ростка.  Сами координатные формы $\Phi$ -- это вещественные полиномы, каждый из которых есть линейная комбинация
вещественных мономов вида
$$ z_1^{\gamma_1} \dots z_p^{\gamma_p} \, \bar{z}_1^{\bar{\gamma}_1} \dots \bar{z}_p^{\bar{\gamma}_p} \,  u_{1}^{\delta_{1}} \dots u_{q}^{\delta_{q}}$$
Используя упомянутые голоморфно полиномиальные преобразования, можно добиться выполнения следующих двух условий:
\begin{eqnarray*}
(I)   \mbox{   координаты всех форм не содержат мономов вида} \quad\quad\quad\quad\quad\quad\qquad\qquad \\
\nonumber
z_1^{\gamma_1} \dots z_p^{\gamma_p} \, u_{1}^{\delta_{1}} \dots  u_{q}^{\delta_{q}}  \mbox{  и  вида } \bar{z}_1^{\gamma_1} \dots \bar{z}_p^{\gamma_p} \, u_{1}^{\delta_{1}} \dots  u_{q}^{\delta_{q}} \quad  \mbox{ ни при каких}  \; \gamma    \mbox{ и } \delta,\\
(II) \quad  \mbox{для любого} \; 1 \leq J \leq p \; \mbox{ни одна из координат формы} \; \Phi_J \quad\quad\quad\quad\qquad\\
\nonumber \mbox{ не содержит слагаемых вида} \qquad c \; \phi(z,\bar{z},u) \, u_{1}^{\delta_{1}} \dots  u_{p}^{\delta_{p}}
 \nonumber   \mbox{   таких что } \qquad \qquad\\
  \phi(z,\bar{z},u)   \; \mbox{это скалярная координата некоторой формы  }   \Phi_j  \mbox{  при }  j <J ,   \mbox{ а } \\
   c \; \mbox{-- это ненулевая константа}. \qquad \qquad\qquad\qquad \qquad\qquad\qquad \qquad\qquad\qquad
\end{eqnarray*}

Отметим, что первое условие -- это условие отсутствия плюригармонических слагаемых.   Второе условие, в отличие от первого,
носит рекуррентный характер. При этом формулировку "не содержит слагаемых" $\;$  следует понимать так, что из мономов, составляющих
координатную форму, нельзя составить выражение указанного вида.

\vspace{3ex}

 Для каждого фиксированного веса $\nu$ обозначим через  $\mathcal{N}_{\nu}$ вещественное линейное пространство
 многочленов веса $\nu$, удовлетворяющих условиям (I) и (II).

\vspace{7ex}

{\bf Определение  (аналитическое):}   Фиксируем $\mu=((\mu_1,n_1), \dots, (\mu_p,n_p))$ -- некоторый вес на переменных группы $z$.  Пусть в окрестности точки $\xi \in M \subset \mathbf{C}^N$ можно ввести координаты так, что локальное уравнение $M$ имеет вид (\ref{1}),  где координатные формы $\Phi$ удовлетворяют условиям (I) и (II),  а также линейно независимы.  Тогда мы говорим, что порождающее многообразие $M$ имеет конечный
$\mu$-тип,  равный $m=((m_1,k_1), \dots, (m_q,k_q))$.

\vspace{3ex}

Отметим, что поскольку минимальный вес переменных группы $z$ -- это $\mu_1$,  а формы $\Phi$ не содержат плюригармонических членов, то минимальный разрешенный вес -- это $2 \, \mu_1$. Таким образом,
$$2 \, \mu_1 \leq m_1 < \dots < m_q.$$

Если в качестве веса выбрать $\mu=((1,n))$, то мы получим модифицированное аналитическое определение Блума-Грэма.

\vspace{3ex}

Прежде чем перейти к геометрическому определению  напомним обычное  геометрическое определение {\it конечности} типа по Блуму-Грэму. Это определение не привязано к выбору каких-либо весов и оно остается у нас без изменений.

\vspace{3ex}

Пусть $D_1$ -- это распределение комплексных касательных, определенное на $M$ в окрестности $\xi$, т.е. $D_1=T^c_M$.  Это распределение можно задать с помощью базисного набора из $2 \, n$ гладких вещественных векторных полей. Далее определим бесконечную последовательность распределений $D_{\nu}$, заданных индуктивно: 
$$D_{\nu+1}=[D_{\nu},D_1]+D_{\nu}, \quad \nu=1,2,...  .$$
 Пусть, далее, $D_{\nu}(\xi)$ - это значение $D_{\nu}$ в точке $\xi$.  Таким образом
$$T^c \; M_{\xi} =D_{1}(\xi) \subset D_{2}(\xi)  \subset  ... \subset  D_{\nu}(\xi) \subset ...$$
Поскольку эта неубывающая последовательность состоит из подпространств $T M_{\xi}$, то она  на каком-то шаге стабилизируется. Если это последнее подпространство совпадает с $T M_{\xi}$, то мы говорим, что $M$ в точке $\xi$ является многообразием {\it конечного типа}, если нет -- {\it бесконечного}.

\vspace{3ex}

Предполагаем далее, что $M$ в точке $\xi$ является многообразием конечного типа.  Теперь возвращаемся в нашу "взвешенную" $\;$ ситуацию.
Ясно, что для многообразия, заданного уравнениями (\ref{1}),  подпространство $D_1(0)$ совпадает с  пространством, порожденным дифференцированиями на пространстве  $\mathbf{C}^n$ координаты $z$ и получает разложение, соответствующее весу $\mu$.
В соответствии с нашим соглашением дифференцирования по координатам получают веса, равные весам координат, но с минусом.
В силу конечности типа в точке описанный индуктивный процесс за некоторое число шагов дает все касательное пространство в точке.
Таким образом, каждый вектор касательного пространства получает выражение через дифференцирования по координатам пространства $\mathbf{C}^n$.  В результате мы получаем
$$T^c_0 M = \mathcal{D}_1(0) \subset \mathcal{D}_2(0) \subset \dots \subset  \mathcal{D}_{max}(0)=T_0 \, M,$$
где подпространство  $ \mathcal{D}_{\nu}(0)$ образовано дифференцированиями веса $-\nu$. 
Пусть $2 \leq m_1 < m_2 < \dots < m_l$ -- это значения тех весов, для которых имел место рост размерности, т.е.
$ \dim \,  \mathcal{D}_{m_j}(0) >   \dim \,  \mathcal{D}_{m_j-1}(0)$, а через $k_j$ обозначим величину этого роста, т.е.
$ \dim \,  \mathcal{D}_{m_j}(0) -   \dim \,  \mathcal{D}_{m_j-1}(0) >0$.

Набор $m=((m_1,k_1), \dots, (m_q,k_q))$ дает второе, {\bf геометрическое определение} $\mu$-типа.  В соответствии с \cite{MS} эти определения эквивалентны. Отметим, что если бы  тип многообразия в точке был бы бесконечным, то как и в стандартном случае,
мы могли бы говорить о $\mu$-типе в точке, дополняя наш набор парой $(\infty,d)$, где
$d=K-(k_1+\dots +k_l)$ -- дефект в точке (см. также \cite{MS2}).  В нашем случае конечного типа дефект равен нулю.

\vspace{3ex}

{\bf Утверждение 1: (см.\cite{MS})}  \\
(a) Геометрическое  и аналитическое определения  $\mu$-типа ростка  эквивалентны.\\
(b)  $\mu$-тип ростка инвариантен относительно локально голоморфных преобразований, сохраняющих весовое разложение комплексной касательной в центре ростка.\\
(c) Если росток имеет конечный тип $m=((m_1,k_1), \dots, (m_q,k_q))$, то для всех $\alpha=1,\dots,p$
$$k_{\alpha} \leq \dim \, \mathcal{N}_{\alpha}.$$

\vspace{3ex}

Отметим, что между распределением весов в переменной $z$ -- величиной $\mu$ -- и распределением весов в переменной $w$ -- величиной $m$ -- имеется существенная разница. Если $\mu$ мы назначаем по собственному произволу, то при фиксированных $M_{\xi}$ и $\mu$ величина $m$, т.е. $\mu$-тип ростка, восстанавливается однозначно. Поэтому можно написать, что $m=m(M_{\xi}, \mu)$.

\vspace{3ex}

{\bf 3. Младшие компоненты отображения}

\vspace{3ex}

Пусть имеется локально обратимое голоморфное отображение  $\chi=(f,g)$ вида
$$(z_{\alpha} \rightarrow f_{\alpha}(z,w), \;\; w_{\beta} \rightarrow g_{\beta}(z,w)) , \; \alpha=1,\dots, p, \;\; \beta = 1,..., q,$$
ростка в начале координат $M_0$ конечного типа $m$, заданного уравнениями в приведенной форме
\begin{eqnarray} \label{M1}
  v_{\beta} = \Phi_{\beta}(z,\bar{z},u)+ F_{\beta}(z,\bar{z},u), \quad j=1,\dots,q,
\end{eqnarray}
в другой такой же росток $\tilde{M}_0$
\begin{eqnarray} \label{M2}
  v_{\beta} = \tilde{\Phi}_{\beta}(z,\bar{z},u)+\tilde{F}_{\beta}(z,\bar{z},u), \quad j=1,\dots,q,
\end{eqnarray}
где $F_{\beta}$  и  $\tilde{F}_{\beta}$ -- это $o(m_{\beta})$,
сохраняющее начало координат и весовое разложение комплексной касательной в нуле.

Рассмотрим  поверхности, заданные уравнениями
\begin{eqnarray} \label{Q}
Q=\{v_{\beta} = \Phi_j(z,\bar{z},u)\},  \quad  \tilde{Q}=\{ v_{\beta} = \tilde{\Phi}_{\beta}(z,\bar{z},u)\}, \quad {\beta}=1,\dots,q.
\end{eqnarray}
Будем говорить, что $Q$ и $\tilde{Q}$ -- {\it $\mu$-модельные поверхности} ростков ${M}_0$ и $\tilde{M}_0$.

\vspace{3ex}

В дальнейшем мы будем использовать разложения
$$f_{\alpha}=\sum \, f_{\alpha \gamma}, \;\;  g_{\beta}=\sum \, g_{ \beta \gamma}, \;\;
F_{\beta}= \sum \, F_{{\beta} \gamma}, \; \; \tilde{F}_{\beta}= \sum \, \tilde{F}_{{\beta}\gamma},$$
где $  f_{\alpha \gamma}, \; g_{\beta \gamma}, \; F_{\beta \gamma}, \; \tilde{ F}_{\beta \gamma}$ - компоненты веса $\gamma$.

В этих обозначениях то, что отображение сохраняет неподвижным начало координат и разложение комплексной касательной на весовые компоненты, означает, что
$$f_{\alpha}=C_{\alpha} \, z_{\alpha} + \tau_{\alpha}(z,w) + o(\mu_{\alpha}), $$
где $C_{\alpha}$ -- невырожденное линейное преобразование пространства $\mathbf{C}^{n_{\alpha}}$, а $\tau_{\alpha}(z,w)$ -- это вектор-полином веса $\mu_{\alpha}$, который может зависеть от $z_{\gamma}$  лишь при $\gamma <\alpha$ и от тех $w_{\delta}$, т.ч. $m_{\delta} < \mu_{\alpha}$. Таким образом, $\tau_{\alpha}(z,w)$  не содержит линейных членов.
Положим
\begin{eqnarray*}
C \, z =(C_1 \, z_1, \dots, C_p \, z_p),   \quad \tau(z,w)=(\tau_1( z,w), \dots, \tau_p(z,w))\\
 \rho \, w =(\rho_1 \, w_1, \dots, \rho_q \, w_q), \quad \theta(z,w)=(\theta_1( z,w), \dots, \theta_q(z,w)) .
\end{eqnarray*}

Записывая, что образ  $M_0$ принадлежит $\tilde{M}_0$, получаем тождество
\begin{eqnarray} \label{MI}
 {\rm Im}\, g = \tilde{\Phi}(f,\bar{f},{\rm Re}\, g)+\tilde{F}(f,\bar{f},{\rm Re}\, g) \;
\mbox{при } w=u+i (\Phi+F).
\end{eqnarray}

\vspace{3ex}

Рассмотрим младшие компоненты тождества (\ref{MI}).

\vspace{5ex}

Начнем с группы переменных $w_1$. В весах от 1 до $(m_1-1)$ получаем  ${\rm Im}\, g_{1\nu} =0$, где $1\leq \nu \leq  m_1-1$. Учитывая, что однородная форма веса $\nu < m_1$ -- это голоморфная форма переменных группы $z$, заключаем, что $ g_{11} =g_{12} =...=g_{1 (m_1-1)}=0$.

 В весе $m_1$ имеем $g_{1 m_1}=a(z)+\rho_1 \, w_1, \quad  f_{\alpha}=C_ {\alpha} \,( z_{\alpha} + \tau_{\alpha}(z,w))$, где $\rho_1$ и $C_{\alpha}$ -- линейны и обратимы, $a(z)$ -- голоморфная однородная форма веса $m_1$. Получаем
$$ {\rm Im}\,(a(z)+\rho_1 \, (u_1+i\,\Phi_1))= \tilde{\Phi}_1(C\,(z +\tau(z,w)),\overline{C\,(z +\tau(z,w))}).$$
Поскольку $m_1$ -- это младший вес в группе весов $w$, то в тех координатах $\tau$, от которых зависит $\Phi_1$,  не могут присутствовать переменные группы $w$.  Т.е. $\tau$ может зависеть от $w$, но не для тех координат, от которых зависит $\Phi_1$.

  Отделяя в этом соотношении компоненту, голоморфную по $z$, и учитывая, что $\Phi_1$ и $\tilde{\Phi}_1$ не содержат голоморфных слагаемых, получаем $a(z)=0$.  Из линейной компоненты по $u_1$  получаем ${\rm Im}\, \rho_1 \, u_1 =0$.  После чего получаем следующее соотношение
\begin{equation}\label{G01}
\rho_1 \Phi_1(z,\bar{z})= \tilde{\Phi}_1(C\,(z +\tau(z)),\overline{C\,(z +\tau(z))})
\end{equation}
Отметим, что это соотношение равносильно тому, что отображение
$(z \rightarrow C \, (z +\tau(z)), \quad w_1 \rightarrow \rho_1 \, w_1)$ переводит "усеченную" $\;$ модельную поверхность $Q(1)=\{v_1=\Phi_1(z,\bar{z})\}$ пространства ${\bf C}^{n+k_1}$, имеющую тип $(m_1,k_1)$, в другую "усеченную" $\;$ поверхность
$\tilde{Q}(1) =\{v_1=\tilde{\Phi}_1(z,\bar{z})\}$ такого же типа.

Переходим к координате $w_2$. Компоненты $g_{2\nu}$, где $\nu < m_2$, -- это выражения вида $\sum \psi_{\alpha\beta}\, (z,w_1)$, где $\psi_{\alpha\beta}\, (z,w_1)$ -- голоморфная полилинейная форма веса $\alpha$ по $z$ и $\beta$ по $w_1$, причем $\alpha +m_1 \,\beta  =\nu$.
Компоненты тождества (\ref{MI})  весов $\nu < m_2$ дают
$${\rm Im}\,(\sum \psi_{\alpha\beta}\, (z,u_1+i\,\Phi_1))=0,$$
Откуда следует, что $g_{2\nu}=0$ при $\nu < m_2$. Это можно было бы доказать непосредственно из полученного тождества. Воспользуемся, однако, другим способом. Действительно, $g_{2\nu}$ -- это голоморфная функция на "усеченном" $\;$ порождающем многообразии конечного типа $Q_1=\{(z,w_1): v_1=\Phi_1(z,\bar{z})\}$, чья мнимая часть равна нулю. Поэтому $g_{2\nu}$  постоянна, а т.к. её вес больше нуля, то это ноль.

В весе $m_2$ имеем
$g_{2 m_2}= \sum \psi_{\gamma \delta}\, (z,w_1)  + \rho_2 \, w_2$, где $\psi_{\gamma\delta}$ имеет вес $\gamma$ по $z$ и степень $\delta$ по $w_1$, $\rho_2$ -- линейна,  причем  $\gamma+ m_1 \,\delta=m_2$,  и при этом имеет место соотношение
\begin{eqnarray*}
{\rm Im}\,\left(\sum \psi_{\gamma \delta}\,(z,u_1+i\,\Phi_1) + \rho_2 (u_2+i\,\Phi_2(z,\bar{z},u_1))\right)=\\
 \tilde{\Phi}_2(C\,(z +\tau(z,w_1)),\overline{C\,(z +\tau(z,w_1))}, \rho_1 \, u_1).
\end{eqnarray*}

Отделяя члены, линейные по $u_2$,  получаем ${\rm Im}( \rho_2) u_2=0$, т.е. линейное отображение $ \rho_2$ -- вещественно. Таким образом, соотношение приобретает вид
\begin{eqnarray*}
{\rm Im}\,\left(\sum \psi_{\gamma \delta}\,(z,u_1+i\,\Phi_1(z,\bar{z}) )\right)=  - \rho_2 \, \Phi_2(z,\bar{z},u_1)\\
+\tilde{\Phi}_2(C\,(z +\tau(z,u_1+i \,\Phi_1)),\overline{C\,(z +\tau(z,u_1+i\,\Phi_1(z,\bar{z})))}, \rho_1 \, u_1).
\end{eqnarray*}
В силу условия (I) его правая часть не содержит слагаемых, голоморфных по $z$. Полагая $\bar{z}=0$, получаем, что если $\gamma \neq 0$, то $\psi_{\gamma \delta}\,(z,u_1)=0$, а
$\psi_{0 \delta}\,(u_1)$ -- вещественная форма веса $m_2$, которую мы переобозначим через
$\rho_2 \, \theta_2(w_1)$, т.е. $g_{2 m_2}=\rho_2 (w_2 + \theta_2(w_1))$. Теперь соотношение приобретает вид
\begin{eqnarray*}
\rho_2 \,\Phi_2(z,\bar{z},u_1)= \qquad\qquad\qquad\qquad\qquad\qquad\qquad\qquad\\
\tilde{\Phi}_2(C\,(z +\tau(z,u_1+i \, \Phi_1)),\overline{C\,(z +\tau(z,u_1+i \, \Phi_1))}, \rho_1 \, u_1) +{\rm Im}\,\theta_2(u_1+i \, \Phi_1)
\end{eqnarray*}

Отметим, что это соотношение равносильно тому, что отображение
$$(z \rightarrow C \, (z  +\tau(z,w)), \quad w_1 \rightarrow \rho_1 \, w_1, \quad w_2 \rightarrow \rho_2 \, (w_2 + \theta_2(w_1)))$$
переводит вторую "усеченную" $\;$ порождающую поверхность
$Q(2) =\{v_1=\Phi_1(z,\bar{z}), \, v_2=\Phi_2(z,\bar{z},u_1)\}$ пространства ${\bf C}^{n+k_1+k_2}$ типа
$((m_1,k_1),(m_2,k_2))$ в другую "усеченную" $\;$ поверхность
$\tilde{Q}(2) =\{v_1=\tilde{\Phi}_1(z,\bar{z}), \, v_2=\tilde{\Phi}_2(z,\bar{z},u_1)\}$ того же типа.

И так далее до последней весовой группы, соответствующей $w_q$. Сформулируем полученный результат.

\vspace{3ex}

{\bf Утверждение 2:} Пусть
$$\chi = (z_{\alpha} \rightarrow f_{\alpha}(z,w), \; w_{\beta} \rightarrow g_{\beta}(z,w) , \quad \alpha=1, \dots, p, \quad {\beta} = 1,..., q)$$
обратимое голоморфное отображение ростка (\ref{M1})  на другой такой росток (\ref{M2}), т.ч. его действие на комплексной касательной в нуле сохраняет её разложение на компоненты весов $\mu$. Тогда\\

(a) Это отображение имеет вид
\begin{eqnarray}  \label{T2}
\nonumber  (z_{\alpha} \rightarrow  C_{\alpha} \, (z_{\alpha}  + \tau_{\alpha}(z,w))+o(\mu_{\alpha}),  \;
w_{\beta} \rightarrow \rho_{\beta} \, (w_{\beta} + \theta_{\beta} (w_1,...,w_{\beta-1})) + o(m_{\beta})  ), \\
\nonumber \mbox{ где } C_{\alpha} \in GL(n_{\alpha},{\bf C}), \; \rho_{\beta} \in GL(k_{\beta},{\bf  R})  \mbox{ и где  } [\tau_{\alpha}]=\mu_{\alpha}, \;\; [\theta_{\beta}]=m_{\beta},\\
   \mbox{ причем для всех } \beta=1,...,q   \qquad\\
\nonumber  \tilde{\Phi}_{\beta}( C\,(z +\tau(z,u+i \, \Phi)),\overline{C\,(z+\tau(z,u+i \, \Phi))},\rho \,( u+{\rm Re} \,(\theta(u+i\,\Phi))))  =\\
\nonumber  \rho_{\beta} \, (\Phi_{\beta}(z,\bar{z},\, u) +{\rm Im}\,\theta_{\beta}(u+i \, \Phi)).
\qquad\qquad\qquad\qquad
\end{eqnarray}
(b) При этом "квазилинейное" $\;$ отображение
\begin{equation}\label{G10}
 (z\rightarrow C \,(z  + \tau(z,w)) ,  \; w_{\nu} \rightarrow \rho_{\nu} \,(w_{\nu} + \theta_{\nu}(w_1,...,w_{\nu-1}))  , \; \nu = 1,...,\beta)
\end{equation}
переводит $\mu$-модельную поверхность $Q$ в $\mu$-модельную поверхность $\tilde{Q}$ (см. (\ref{Q})).

\vspace{3ex}
Более того, для каждого $\beta=1,...,q$ усеченное отображение
$$ (z\rightarrow C \,z + \tau(z,w),  \; w_{\nu} \rightarrow \rho_{\nu} \,w_{\nu} + \theta_{\nu}(w_1,...,w_{\nu-1})  , \; \nu = 1,...,\beta)$$
переводит усеченную модельную поверхность
$$Q(\beta)=\{ v_{\nu}=\Phi_{\nu}, \; \nu =1,...,\beta\}$$
пространства ${\bf C}^{n+k_1+...+k_{\beta}}$ в соответствующую усеченную модельную поверхность.
$$\tilde{Q}(\beta)=\{ v_{\nu}=\tilde{\Phi}_{\nu}, \; \nu =1,...,\beta \}.$$
Пункт (b) теоремы это частный случай этого утверждения.

\vspace{3ex}

Если положить $\tilde{Q}=Q$,  то отображения (\ref{G10}) c условием (\ref{T2}) образуют некоторую подгруппу $G_0$ автоморфизмов $Q$, сохраняющих начало координат и весовое разложение комплексной касательной в нуле.
Это автоморфизмы $Q$ вида
\begin{eqnarray}\label{G0}
\nonumber
z_{\alpha} \rightarrow C_{\alpha} \, (z_{\alpha} + \tau_{\alpha}(z,w)) , \\
w_{\beta} \rightarrow \rho_{\beta} \,(w_{\beta} + \theta_{\beta}(w_1,...,w_{\beta-1})).
 \end{eqnarray}
 Т.е. это автоморфизмы, т.ч. каждая координата сохраняет свой вес.
 Если привлечь введенное нами в пункте 4 понятие составляющей, то $G_0$ можно охарактеризовать так:
это автоморфизмы $Q$, чье разложение на составляющие содержит лишь 0-составляющую.
Именно это фиксирует (\ref{G0}). То, что такое отображение является автоморфизмом $Q$, -- это условия (\ref{T2})
(утверждение 2) при $\tilde{\Phi}={\Phi}$.

\vspace{3ex}
Если рассмотреть совокупность преобразований вида (\ref{G0}) без условий (\ref{T2}), то они образуют
подгруппу полиномиальных автоморфизмов пространства
$$\mathbf{C}^N=\mathbf{C}^{n_1}+\dots+\mathbf{C}^{n_{p}} +\mathbf{C}^{k_1}+\dots+\mathbf{C}^{k_{q}},$$
Эта подгруппа $\mathcal{G}_0$ является полупрямым произведением подгруппы треугольных преобразований вида
\begin{eqnarray*}
z_{\alpha} \rightarrow  z_{\alpha} + \tau_{\alpha}(z,w) ,  \quad w_{\beta} \rightarrow  w_{\beta} + \theta_{\beta}(w_1,...,w_{\beta-1}).
 \end{eqnarray*}
и линейных вида
\begin{eqnarray*}
z_{\alpha} \rightarrow C_{\alpha} \, z_{\alpha}, \quad w_{\beta} \rightarrow \rho_{\beta} \,w_{\beta}.
 \end{eqnarray*}

\vspace{2ex}

В стандартной, невзвешенной версии \cite{VB20}, имеется утверждение (теорема 5 п.(f)) о том, что элемент $G_0$ однозначно определяется своим действием на координату $z$. Вот аналог этого утверждения.

\vspace{3ex}

{\bf Утверждение 3:} Если модельная поверхность $Q$ имеет в нуле конечный $\mu$-тип
и имеется автоморфизм $ (f(z,w),g(z,w)) \in G_0$, т.ч. $f(z,w)=z$, то $g(z,w)=w$, т.е. это тождественное отображение.\\
{\it Доказательство:}  Положим $\tau(z,w)=0, \; C \, z = z$ и запишем соотношение (\ref{T2}) для $\beta=1$, получим
$\rho_1 \Phi_1(z,\bar{z})= {\Phi}_1(z,\overline{z})$. Откуда из линейной независимости координат $\Phi_1$ следует,
что $\rho_1 \, w_1=w_1$.  Запишем соотношение (\ref{T2}) для $\beta=2$, получим
$$   {\Phi}_2(z,\overline{z}, u_1)= \rho_2 \,\Phi_2(z,\bar{z},u_1)+{\rm Im}\,\theta_2(u_1+i \, \Phi_1)$$
Из того, что $\Phi_2$ записана в приведенной форме, получаем $\theta_2=0$ и $\rho_2 \, w_2=w_2$.  И так далее до
$\beta=q$. Утверждение доказано.

\vspace{3ex}

Это утверждение можно проинтерпретировать как наличие некоторой параметризации группы $G_0$.
Всякий элемент $\chi \in G_0$ однозначно определяется набором параметров $(C,\tau,\rho,\theta)$, где набор параметров
связан алгебраическими соотношениями (\ref{T2}). Утверждение 3 означает, что параметры $(\rho,\theta)$ однозначно
определяются параметрами $(C,\tau)$, т.е. $\rho=\rho(C,\tau), \;\theta=\theta(C,\tau)$.  Алгебраическое подмножество в пространстве $(C,\tau,\rho,\theta)$, определяемое соотношениями (\ref{T2}) имеет однозначную проекцию в некоторое алгебраическое подмножество $\mathcal{C}$ пространства $(C,\tau)$.  В результате $\mathcal{C}$ получает структуру алгебраической группы, действующей в $\mathbf{C}^N$ и изоморфной $G_0$.

\vspace{3ex}

{\bf 4. Взвешенная конструкция Пуанкаре}

\vspace{3ex}

Введем в линейном пространстве $\mathcal{V}$ наборов голоморфных ростков в окрестности нуля вида $\chi=(f_1,\dots,f_p;g_1,\dots,g_q)$ прямое разложение на {\it составляющие}  $\mathcal{V}=\sum \,  \mathcal{V}_{\nu}$,
где $\mathcal{V}_{\nu}$ состоит наборов следующих весов $((\mu_1+\nu,\dots,\mu_p+\nu),(m_1+\nu, \dots, m_q+\nu))$.
Соответственно мы можем написать $\chi=\sum \, \chi^{(\nu)}=(f,g)=\sum \, (f^{(\nu)},g^{(\nu)})$,  при этом
$$ f^{(\nu)}=(f_{1 (\mu_1+\nu)}, \dots, f_{p (\mu_p+\nu)}), \quad   g^{(\nu)}=(g_{1 (m_1+\nu)}, \dots, g_{q (m_q+\nu)}).$$

Пусть младшие (модельные) члены уравнений (\ref{M1}) и (\ref{M2}) совпадают, т.е. $\tilde{\Phi}=\Phi $ и $\tilde{Q}=Q$.
Как было показано выше, младшая составляющая отображения $\chi=\sum \, \chi^{(\nu)}$ поверхности  (\ref{M1}) на  (\ref{M2}) при условии сохранения начала координат и весового разложения комплексной касательной -- это $\chi_{0}$, т.е. набор координатных функций имеет веса вида $((\mu_1,\dots,\mu_p),(m_1, \dots, m_q))$.
При этом $\chi_{0}$ -- это элемент линейной алгебраической группы $G_0$ и он, как было показано, задается своей системой параметров $\lambda_0=(C,\tau,\rho,\theta)$, т.е. $\chi_{0}=\chi_{0}(\lambda_0)$.
 Ясно, что любое такое отображение $\chi=\chi_0+\chi_1+\dots$ может быть представлено в виде композиции $ \varphi \circ \psi$, где $\psi=\varphi^{-1} \circ \chi$ -- это отображение, у которого 0-составляющая -- это тождественное отображение, т.е. $\psi = Id +\psi^{(1)}+\psi^{(2)}+\dots$.

Применим конструкцию Пуанкаре, чтобы дать оценку размерности семейства отображений вида $\psi = Id +\psi^{(1)}+\psi^{(2)}+\dots$. Пространство $\mathcal{V}_1+\mathcal{V}_2+\dots$ обозначим через $\mathbf{V}_{+}$.
Рассмотрим теперь тождество (\ref{MI}) и выделим в нем $\nu$-ю составляющую. Получаем соотношение
\begin{eqnarray} \label{GOP}
\nonumber
 - {\rm Im} \, g^{ (\nu)} + d \, \Phi(z,\bar{z},u)\, (f^{(\nu)},\bar{f}^{(\nu)}, {\rm Re}\, g^{ (\nu)}) = \qquad \qquad\qquad\qquad\\
 \mbox{членам, зависящим от составляющих  }  f^{(\iota)},  \, g^{(\iota)} \mbox{  при }  \iota < \nu,\\
\nonumber \mbox{  где  } w=u+i \, \Phi(z,\bar{z},u).
\end{eqnarray}
Пусть $\mathcal{K}$ -- ядро линейного оператора
\begin{equation}\label{L}
  \mathcal{L}(f,g)= - {\rm Im} \, g + d \, \Phi(z,\bar{z},u)\, (f,\bar{f}, {\rm Re}\, g),  \mbox{  где  } w=u+i \, \Phi(z,\bar{z},u),
\end{equation}
действующего на $\mathbf{V}_{+}$.  Это линейное подпространство $\mathbf{V}_{+}$, которое можно разложить по составляющим $\mathcal{K}^{(1)}+\mathcal{K}^{(2)}+\dots$.  Независимо от конечномерности ядра $\mathcal{K}$ отдельно каждая его составляющая конечномерна, т.к. ее координатные проекции являются подпространствами полиномов фиксированного веса.

Основное наблюдение, на котором основаны применения конструкции Пуанкаре в $CR$-геометрии, состоит в следующем. Условие, что векторное поле  в окрестности начала координат
 $$X=2 \,{\rm Re} \left(  f \, \frac{\partial}{\partial z} + g \, \frac{\partial}{\partial w} \right)$$
принадлежит алгебре Ли ${\rm aut} \, Q$  инфинитезимальных голоморфных автоморфизмов модельной поверхности $Q$, -- это $\mathcal{L}(f,g)=0$.

Рассмотрим соотношение (\ref{GOP}) как рекуррентное соотношение для вычисления последовательных составляющих
отображения $\psi$.

Первый шаг -- выделяем 1-составляющую (\ref{GOP}). Она имеет вид
$\mathcal{L}(f^{(1)},g^{(1)})=T_1(\lambda_0)$. Мы видим, что для однозначного определения $(f^{(1)},g^{(1)})$ достаточно зафиксировать элемент  $\lambda_1 \in \mathcal{K}^{(1)}$, и мы можем написать $(f^{(1)},g^{(1)})=(f^{(1)}(\lambda_1),g^{(1)}(\lambda_1))$. Однако для разрешимости полученной неоднородной линейной системы следует добавить условие разрешимости (условие попадания $T_1(\lambda_0)$ в образ $\mathcal{L}(f^{(1)},g^{(1)})$. Это условие $C_1(\lambda_0 ,\lambda_1)=0$ есть вещественно алгебраическое соотношение между $\lambda_0 \in G_0$ и  $\lambda_1 \in \mathcal{K}^{(1)}$

Второй шаг -- выделяем 2-составляющую (\ref{GOP}). Она имеет вид
$\mathcal{L}(f^{(2)},g^{(2)})=T_2(\lambda_1)$, где $T_2$ -- вещественный вектор-полином от $\lambda_1$. Мы видим, что при фиксированном $\lambda_1$ для однозначного определения $(f^{(2)},g^{(2)})$ достаточно зафиксировать элемент  $\lambda_2 \in \mathcal{K}^{(2)}$. А также возникает новое условие разрешимости
$C_2(\lambda_0 ,\lambda_1,\lambda_2)=0$. И так далее до бесконечности.

Если ${F}(z,\bar{z},u)=\tilde{F}(z,\bar{z},u)=0$, т.е. отображение -- это автоморфизм модельной поверхности, то
все условия согласования исчезают, т.е. $C_{\nu}=0$ для всех $\nu$. Таким образом, система параметров, задающая автоморфизм $Q$ вида $\chi=\chi_0+\chi_1+\dots$, совпадает с $G_0 \cup \mathcal{K}$.

В общем случае мы получаем следующее описание системы параметров
\begin{equation}\label{L}
\Lambda(M_0) = \{(\lambda_0,\lambda_1,\lambda_2,\dots) \in G_0 \cup \mathcal{K}: C_{\nu}(\lambda_0,\lambda_1,\dots.\lambda_{\nu})=0, \; \nu=1,2, \dots\}.
\end{equation}

Если ядро $\mathcal{K}$ оказалось конечномерным и, тем самым, конечноградуированным, т.е. $\mathcal{K}_{\nu}=0$
при $\nu >d$, то $\Lambda$, как мы видим, есть вещественно алгебраическое подмножество конечномерного вещественного пространства.

Если $\tilde{M}_0=M_0$, т.е. семейство отображений, которое параметризует совокупность параметров   $\Lambda$ -- это автоморфизмы $M_0$, то множество $\Lambda$ естественно отождествляется с некоторой подгруппой автоморфизмов $M_0$. Эту подгруппу, состоящую из автоморфизмов вида $\chi=\chi_0+\chi_1+\dots$, мы обозначим через ${\rm Aut}^{\mu}_{0} \, M_{0}$.   Это в точности те автоморфизмы $M_{0}$, которые сохраняют начало координат и весовое разложение комплексной касательной в нуле. Эту подгруппу стабилизатора нуля будем называть {\it $\mu$-стабилизатором}. Используя соответствие $\lambda \rightarrow \chi(\lambda)$, мы можем написать
$$\Lambda \approx {\rm Aut}^{\mu}_{0} \, M_{0}.$$
Поскольку группа ${\rm Aut}^{\mu}_{0} \, M_{0}$ транзитивно действует на себе левыми сдвигами, то это действие
индуцирует транзитивное действие на $\Lambda$. Если $\mathcal{K}$  конечномерно, то $\Lambda$ -- вещественно алгебраическое многообразие без особенностей.

Соответствующую алгебру Ли обозначим как ${\rm aut}^{\mu}_{0} \, M_{0}$.   Это некоторая подалгебра в полном стабилизаторе ${\rm aut}_{0} \, M_{0}$, который, в свою очередь, есть подалгебра в полной алгебре ${\rm aut} \, M_{0}$.

Подгруппу автоморфизмов $M_0$, состоящую из автоморфизмов вида   $\chi=Id+\chi_1+\dots$, обозначим
$G_{+} ( M_{0})$. Подгруппу автоморфизмов $M_0$, состоящую из автоморфизмов вида   $\chi=\chi_0$ обозначим
$G_{0} ( M_{0})$. 

 \vspace{3ex}

Перед  общей формулировкой полученного результата приведем частные утверждения, имеющее отношение к модельной поверхности $Q$.

\vspace{3ex}

{\bf Утверждение 4:}\\
(a) $G_{+}( Q_{0})\approx \mathcal{K}$, причем $\mathcal{K}$ -- это линейное пространство.\\
(b) $ {\rm Aut}^{\mu}_{0} \, Q_{0} \approx \Lambda(Q_0) = G_0 \sqcup \mathcal{K}.$\\
(c) $ {\rm Aut}^{\mu}_{0} \, Q_{0} = G_0 ( M_{0}) \ltimes  G_{+}( Q_{0})$ -- полупрямое произведение.

\vspace{3ex}

Итак,\\
-- пусть $\mu$ -- произвольная градуировка переменной $z$,\\
-- пусть имеется два ростка в начале координат $M_0$ и $\tilde{M}_0$ конечного типа $m$, заданного уравнениями в приведенной форме
\begin{eqnarray*}
  M_0 =\{v_{\beta} = \Phi_{\beta}(z,\bar{z},u)+ F_{\beta}(z,\bar{z},u)\}, \\
\tilde{M}_0=\{v_{\beta} = \Phi_{\beta}(z,\bar{z},u)+\tilde{F}_{\beta}(z,\bar{z},u)\},\\
 \quad j=1,\dots,q,  \quad   F_{\beta}=o(m_{\beta}), \;\;  \tilde{F}_{\beta}=o(m_{\beta}).
\end{eqnarray*}
-- рассмотрим

Наши рассуждения доказывают следующую теорему.

\vspace{3ex}

 {\bf Теорема 5:}   Пусть $\mu$ -- произвольная градуировка переменной $z$, пусть $Q_0$ -- $\mu$-модельная поверхность ростка $M_0$, тогда\\
(a) $G_{0} ( M_{0})$  -- это подгруппа в $G_{0} ( Q_{0})$.\\
(b) Множество $\Lambda( M_{0})$, параметризующее семейство обратимых отображений
$M_0$ на $\tilde{M}_0$ вида $\chi = \chi_0 +\chi^{(1)}+\chi^{(2)}+\dots$ имеет вид (\ref{L}).\\
(c) Если ${\rm aut}^{\mu}_{0} \, Q_{0} < \infty$, то $\Lambda( M_{0})$ -- это неособое вещественно алгебраическое множество и $\dim \Lambda( M_{0})  \leq  \dim \,{\rm aut}^{\mu}_{0} \, Q_{0}.$\\
(d) $\dim \,{\rm aut}^{\mu}_{0} \, M_{0} \leq  \dim \,{\rm aut}^{\mu}_{0} \, Q_{0}.$\\

 \vspace{3ex}

 В связи с этой теоремой уместно задать следующий вопрос. Когда  $\dim \,{\rm aut}^{\mu}_{0} \, Q_{0}$ конечна?

 \vspace{3ex}

 {\bf Утверждение 6:}     $\dim \,{\rm aut}^{\mu}_{0} \, Q_{0} < \infty$ тогда и только тогда, когда
 $\mu$-модельная поверхность $Q$ имеет конечный тип и голоморфно невырождена.\\
 {\it Доказательство:}  Достаточность следует из известной теоремы (см. \cite{BER}). Покажем необходимость.
 Необходимость голоморфной невырожденности -- очевидна.  Если $\mu$-модельная поверхность $Q$ имеет бесконечный тип, то это означает, что среди координатных форм какого-то веса имеется линейная зависимость. После линейного преобразования соответствующей группы переменных $w$ среди координатных форм появятся тождественные нули.
 Это сразу дает бесконечномерность. Утверждение доказано.

\vspace{3ex}
Это утверждение, как и его доказательство вполне аналогичны тому, которые имеются при обычном невзвешенном подходе \cite{VB20}. Отметим, что условия утверждения 6 -- это также и критерий конечномерности полной алгебры  ${\rm aut}\, Q_{0}$.

 \vspace{3ex}

 В связи с этим утверждением мы напоминаем  следующее определение.  Росток называется {\it невырожденным},
 если он имеет конечный тип и голоморфно невырожден.  Это определение не зависит от выбора $\mu$. Поэтому невырожденность по отношению к одному весу означает невырожденность по отношению ко всем весам.  Невырожденность $\mu$-модельной поверхности $Q$ влечет невырожденность соответствующего ростка. Обратное  неверно.

 \vspace{3ex}

 {\bf 5. Взвешенная модельная поверхность и её автоморфизмы}

\vspace{3ex}
Введенные для переменных групп $z$ и $w$ веса естественно продолжаются и на дифференцирования. В результате алгебра Ли всех инфинитезимальных голоморфных автоморфизмов любого ростка становится градуированной алгеброй Ли. Эту градуированную алгебру мы обозначим ${\rm aut}^{\mu} \, Q_{0}$. Введение в обозначение для алгебры автоморфизмов верхнего индекса $\mu$ подчеркивает то, что различные градуировки переменных группы $z$ превращают одну и ту же алгебру автоморфизмов в разные {\it градуированные} алгебры Ли.

Переменным  группы $w_q$ максимального веса $m_q$ соответствуют дифференцирования старшего отрицательного веса $-m_q$. Поэтому можно написать, что $$ {\rm aut}^{\mu} \, Q_{0}= \sum_{\nu=-m_q}^{\infty} g_{\nu}.$$
Мы можем рассмотреть три подалгебры, в сумму которых раскладывается полная алгебра: $g_{-}=\sum_{\nu<0}\, g_{\nu}, \;\; g_0, \;\;  g_{+}=\sum_{\nu>0}\, g_{\nu},$
$$    {\rm aut}^{\mu} \, Q_{0}=g_{-} + g_0 + g_{+}.$$

Непосредственно проверяем справедливость следующего утверждения.

\vspace{3ex}

{\bf Утверждение 7 :}\\
(a)  $g_0$ -- это алгебра Ли группы $G_0$.\\
(b)   Подгруппа $G_0$ содержит 1-параметрическую (градуирующую) подгруппу
\begin{eqnarray} \label{DL}
(z_{\alpha} \rightarrow t^{\mu_{\alpha}} \, z_{\alpha}, \; \; w_{\beta} \rightarrow t^{m_{\beta}} \; w_{\beta}), \quad  t \in {\bf R}^{*}
\end{eqnarray}
Этой подгруппе соответствует векторное поле веса ноль
\begin{eqnarray} \label{GVF}
X_0=2 \,{\rm Re}\,(\sum \, \mu_{\alpha} \, z_{\alpha} \,\frac{\partial}{\partial \,z_{\alpha}}+
\sum m_{\beta} \,w_{\beta} \,\frac{\partial}{\partial \,w_{\beta}}),
\end{eqnarray}
(c)  Если $X =\sum_{\nu=-m_q}^{\infty} X_{\nu} \in  {\rm aut}^{\mu} \, Q_{0},$ то $X_{\nu} \in  {\rm aut}^{\mu} \, Q_{0}$ для каждого $\nu$.\\
(d)  $\mu$-стабилизатор нуля -- это полупрямое произведение $G_0 \rtimes G_{+}$, соответственно
 ${\rm aut}^{\mu}_0 \, Q_0 = g_0+g_{+}$.\\
(e) ${\rm aut}^{\mu} \, Q_{0}$ -- конечномерна тогда и только тогда, когда ${\rm aut}^{\mu} \, Q_{0}$ -- конечноградуирована,
т.е. лишь конечное число весовых компонент $g_{\nu}$ отлично от нуля.

 \vspace{3ex}

 Таким образом, если $Q$  невырождена, то
$$ {\rm aut}^{\mu} \, Q_{0}= \sum_{\nu=-m_q}^{d} g_{\nu},$$
где $d$ -- старшая неотрицательная ненулевая компонента.

\vspace{3ex}

Пусть $\xi =(z=a, \; w=b)$ -- точка $Q$. Рассмотрим первую группу уравнений $Q$, а именно $v_1=\Phi_1(z,\bar{z})$, сделаем
замену $z \rightarrow a+z, \;\; w_1 \rightarrow (b_1+i \,\Phi_1(a,\bar{a}))  +w_1$ и запишем получившиеся уравнения:
$$ v_1 =  \Phi_1(z+a,\bar{z}+\bar{a})-\Phi_1(a,\bar{a})=\Phi_1(z,\bar{z}) + d \,\Phi_1(z,\bar{z})(a,\bar{a})+\dots .$$
Если присвоить параметрам $a$ соответствующие $\mu$-веса, то все слагаемые будут иметь вес $m_1$. Если же считать веса только по переменным $z$, то первое слагаемое имеет вес $m_1$, а все остальные -- строго меньше.  Пользуясь известной процедурой, мы можем эту "усеченную" $\;$  поверхность $Q(1)$ записать в приведенном виде. Для этого первого шага процедура направлена только на то, чтобы убрать плюригармонические члены. Если при этом останутся слагаемые веса меньше чем $m_1$, то мы можем констатировать, что $\mu$-тип изменился и ввести новые веса для переменных группы $w_1$. Если же этого не произошло, т.е. после приведения уравнение $Q(1)$ вернулось к прежнему виду,
то это означает, что существует голоморфный автоморфизм  $Q(1)$ вида
$$z  \rightarrow a+z, \;\; w_1 \rightarrow b_1 +w_1 + P_1(z,a, \bar{a}),$$
переводящий $(z=0,\; w_1=0)$  в $(z=a, \; w_1=b_1+i \, \Phi_1(a,\bar{a}))$, причем вес полинома $P_1$ строго меньше $m_1$.

Далее мы переходим к рассмотрению второй "усеченной" $\;$ поверхности $Q(2)$, заданной уравнениями $Q(1)$  и группой уравнений вида $v_2=\Phi_2(z,\bar{z},u_1)$ в точке  $(z=a, \; w_1=b_1+i \, \Phi_1(a,\bar{a}), \; w_2=b_2+i \, \Phi_2(a,\bar{a},b_1))$.   При этом уравнения $Q(1)$ -- это новые приведенные уравнения. Теперь опять повторяется процедура построения приведенных уравнений $Q(2)$ и вычисления $\mu$-типа в новой точке. Если тип не изменился, то мы получаем полиномиально-треугольный автоморфизм $Q(2)$
вида
\begin{eqnarray*}
z  \rightarrow a+z, \;\; w_1 \rightarrow (b_1+i \, \Phi_1(a,\bar{a})) +w_1 + P_1(z,a,\bar{a}) ,\\
 w_2 \rightarrow (b_2+i \, \Phi_2(a,\bar{a},b_1)) +w_2 + P_2(z,w_1,a, \bar{a},b_1),
\end{eqnarray*}
т.ч. он переводит $(z=0,\; w_1=0,\;w_2=0)$  в $(z=a, \; w_1=b_1+i \, \Phi_1(a,\bar{a}),\;w_2=b_2+i \, \Phi_2(a,\bar{a},b_1))$, причем вес $P_2$ строго меньше $m_2$.
И так далее до последней группы переменных. Построенный в результате этого процесса однозначно определенный полиномиально-треугольный "сдвиг" $\;$ обозначим через $S_{\xi}$.

\vspace{3ex}

Пусть $Q$ -- взвешенная модельная поверхность, чей $\mu$-тип в начале координат равен $m$.  Пусть $Q^m$ -- это совокупность точек $Q$, т.ч. в этих точках $\mu$-тип равен $m$. Итак, для каждой точки
$\xi \in Q^m$ существует единственный автоморфизм $S_{\xi}$, построенный в результате приведения уравнений и который переводит начало координат в $\xi$.  Совокупность этих $"$сдвигов$"$ -- это подгруппа, которую мы обозначим через $GS$, а ее алгебру -- через $gs$.

Чтобы явно описать поля из $gs$, надо продифференцировать полученные треугольно-полиномиальные замены по параметрам (выделить линейные по параметрам слагаемые).

\vspace{3ex}

Пусть $St$ -- стабилизатор начала координат в группе $ {\rm Aut} \, Q_{0}$ голоморфных автоморфизмов $Q$,  а $st$ -- соответствующая алгебра Ли.  Ясно, что $St$ -- это подгруппа, а $st$ -- это подалгебра.  При стандартном (невзвешенном) подходе $st$ -- это сумма всех неотрицательных компонент алгебры (т.е. $g_0+g_{-}$).  При взвешенном подходе это, вообще говоря, не так.  В связи с этим обозначим через $st_{-}$ подалгебру, состоящую из  полей, входящих в $g_{-}$, которые обращаются в ноль в начале координат.  Поля
$$X=(f,g)=(f_1, \dots,f_p; \, g_1,\dots,g_q),$$
принадлежащие $st_{-}$, -- это поля, которые кроме условия касания $Q$
$$
{\rm Im} \, g =  d \, \Phi(z,\bar{z},u)\, (f,\bar{f}, {\rm Re} g),  \mbox{  где  } w=u+i \, \Phi(z,\bar{z},u),
$$
удовлетворяют следующим условиям
\begin{equation} \label{st}
f(0,0)=0, \; g(0,0)=0, \; \mbox{вес } f_{\alpha} < \mu_{\alpha},  \; \mbox{вес } g_{\beta} < m_{\beta}.
\end{equation}
С другой стороны, алгебра $gs$ -- это поля c условием касания, такие что
\begin{equation} \label{gs}
f(z,w)=a=const, \; g_1=b_1+p_1(z,a), \; g_2=b_2+p_2(z,w_1,a,b_1), \dots,
\end{equation}
где полиномы $p_{\beta}$ однозначно восстанавливаются из условия касания, они имеют вес меньший, чем $m_{\beta}$, и обращаются в ноль при $a=0, \; b=0$.

\vspace{1ex}

Итак, получаем
\vspace{1ex}

{\bf Утверждение 8:}\\
(a)   Имеет место  разложение подалгебры $g_{-} = gs +st_{-}$.\\
(b)  $gs$ -- это подалгебра, соответствующая подгруппе сдвигов $GS$.\\
(b)   Орбита начала координат по отношению к полной группе автоморфизмов $ {\rm Aut} \, Q_{0}$ совпадает
с орбитой начала координат по отношению к $GS$.

\vspace{3ex}

Пусть $G_{+}$ -- подгруппа группы автоморфизмов $Q$, состоящая из автоморфизмов, чьи разложения по составляющим имеют вид
$$\chi = Id + \chi^{(1)}+\chi^{(2)} + \dots.$$
Это те замены координат, которые использовались нами при описании конструкции Пуанкаре.
Аналогично через $G_{-}$ обозначим подгруппу группы автоморфизмов $Q$, состоящую из автоморфизмов, чьи разложения по составляющим имеют вид
$$\chi = Id + \chi^{(-1)}+\chi^{(-2)} + \dots+\chi^{(-m_q)}.$$
Если в разложении элементов $G_{+}$ могут, вообще говоря, присутствовать составляющие со сколь угодно большим номером, то разложение $G_{-}$ ограничено снизу числом $-m_q$, где $m_q$ -- старший вес в группе переменных $w$.
Отсюда автоматически следует, что преобразования из $G_{-}$ -- полиномиальны.  Положим $St_{-}=G_{-} \cap St$.

 Имеем следующее утверждение:
\vspace{1ex}

 {\bf Утверждение 9:} (a)  Подгруппе автоморфизмов $G_{+}$ соответствует  алгебра Ли $g_{+}$.\\
 (b)  Подгруппе автоморфизмов $G_{-}$ соответствует  алгебра Ли $g_{-}$.\\
 (c)  Подгруппе автоморфизмов $St_{-}$ соответствует  алгебра Ли $st_{-}$.

\vspace{1ex}

Итак, мы имеем разложение ${\rm aut}^{\mu} \, Q_{0}=gs + st_{-} + g_0 + g_{+}$.
Этим четырем подалгебрам соответствуют четыре подгруппы
$$GS, \quad St_{-}, \quad  G_0, \quad G_{+},$$  которые были описаны выше.

\vspace{1ex}

Относительно подгруппы $St_{-}$ имеет место утверждение, вполне аналогичное утверждению 3.

\vspace{1ex}

{\bf Утверждение 10:} Если модельная поверхность $Q$ имеет в нуле конечный $\mu$-тип
и имеется автоморфизм $ (f(z,w),g(z,w)) \in St_{-}$, т.ч. $f(z,w)=z$, то $g(z,w)=w$, т.е. это тождественное отображение.

Доказательство без изменений переносится из утверждения 3.

\vspace{1ex}

Выше (теорема 5) с помощью рекуррентной конструкции Пуанкаре мы показали, что  размерность  $\mu$-стабилизатора точки $\xi$ в группе автоморфизмов ростка оценивается через размерность  $\mu$-стабилизатора нуля в группе автоморфизмов ее модельной поверхности. Однако, чтобы получить оценку всего стабилизатора необходимо оценить размерность $St_{-}(M_{\xi})$. Для оценки размерности $st_{-}$ имеет место полный аналог теоремы 5.
Однако в схеме доказательства имеется отличие. В то время как доказательство теоремы 5 -- это применение конструкции Пуанкаре к семейству отображений, доказательство теоремы 11 -- это конструкция Пуанкаре, примененная к векторным полям. Причем использование конструкции Пуанкаре для оценки размерностей не семейства отображений, а пространств векторных полей позволяет оценить размерность алгебры возмущенной поверхности через соответствующую размерность для модельной не только для $st_{-}(M_{\xi})$, но и для $g_{-}(M_{\xi})$ и для всей ${\rm aut} \, (M_{\xi})$.

\vspace{3ex}

{\bf Теорема 11:}   Пусть $\mu$ -- произвольная градуировка переменной $z$, пусть $Q_0(\mu)$ -- $\mu$-модельная поверхность ростка $M_0$, тогда
$$\dim \,{\rm aut} \, M_{0} \leq  \dim \,{\rm aut} \, Q_0(\mu).$$
{\it Доказательство:} Пусть $\{v=\Phi(z,\bar{z},u)+\Psi(z,\bar{z},u)\}$ -- росток,  $\{v=\Phi(z,\bar{z},u)\}$ -- модельная поверхность, а $\Psi(z,\bar{z},u)$ -- возмущение.  Т.е. если $\Phi_j$ и $\Psi_j$ -- координаты $\Phi$ и $\Psi$, соответствующие $w_j$, то  $\Phi_j$ имеет вес $m_j$, а $\Psi_j$ есть сумма компонент большего веса  $\Psi_j=\Psi_{j (m_j+1)}+\dots$.
Можно записать уравнение возмущенной поверхности в виде разложения по составляющим, а именно
$v=\Phi + \Psi^{(1)}+\Psi^{(2)}+\dots$, где $\Psi^{(\nu)}=(\Psi_{1 (m_1+\nu)},\dots,\Psi_{q (m_q+\nu)})$. Ясно, что 0-составляющая правой части уравнения -- это $\Phi$.
Пусть
$$X=2\, {\rm Re} \left(\sum f_{i}(z,w) \frac{\partial}{\partial z_{\nu}} +  \sum g_{j}(z,w) \frac{\partial}{\partial w_{j}} \right)
=(f(z,w),g(z,w))  -- $$
поле из $ {\rm aut}\, M_0$ и $X=\sum \, X^{(\nu)}$ -- его разложение по весовым компонентам. Тогда
$$X^{(\nu)}=(f^{(\nu)},g^{(\nu)})=(f_{1(\mu_1+\nu)},\dots,f_{q(\mu_q+\nu)};g_{1(m_1+\nu)},\dots,g_{q(m_q+\nu)})$$
-- набор коэффициентов $X^{(\nu)}$ ($\nu$-составляющая поля).  Записывая условие касания, получаем
\begin{eqnarray} \label{TC}
\nonumber
 \mathcal{L}(f,g;\Phi,\Psi) = -{\rm Im} \, g(z,w) +  2 \,{\rm Re} \,(\partial_z (\Phi(z,\bar{z},u) + \Psi(z,\bar{z},u))(f(z,w)) ) + \\
 \nonumber \partial_u (\Phi(z,\bar{z},u) + \Psi(z,\bar{z},u))({\rm Re} \, g(z,w))=0,\\
 \mbox{где  } w=u+i\,(\Phi(z,\bar{z},u)+\Psi(z,\bar{z},u)).
 \end{eqnarray}
 Это линейное выражение  (оператор) $\mathcal{L}(f,g;\Phi,\Psi)$ есть сумма
 $$\mathcal{L}(f,g;\Phi,\Psi)=L(f,g;\Phi)+L'(f,g;\Phi,\Psi),$$
где  $L(f,g;\Phi)=\mathcal{L}(f,g;\Phi,0)$. Т.е. $L(f,g;\Phi)=0$ -- это запись того, что поле $(f,g)$ касается модельной поверхности $Q$, а все члены  $\mathcal{L}$, зависящие от $\Psi$, собраны в $L'$. Будем записывать соотношение (\ref{TC}) как равенство нулю его составляющих. Ясно, что если показатель $\nu$ достаточно мал ($\nu < -m_q$), то $X^{(\nu)}=0$ и соотношение в такой составляющей не даст соотношений на коэффициенты поля.
Первое содержательное соотношение мы получим для $\nu=-m_q$.
Имеем $X^{(-m_q)}=(0,\dots,0;0,\dots,0,g_{q 0})$, где $[g_{q 0}]=0$, т.е. $g_{q 0}$ -- константа. Нетривиальное соотношение
для $(-m_q)$-составляющей (\ref{TC}) -- это $ {\rm Im} \, g_{q 0}=0$. В следующей $(-m_q+1)$-й составляющей участвуют
$X^{(-m_q+1)}=(f^{(-m_q+1)},g^{(-m_q+1)})$ и $X^{(-m_q)}$. И т.д.,  причем  $\nu$-я составляющая (\ref{TC}) имеет вид
$$L(f^{(\nu)},g^{(\nu)})+ \mbox{члены, зависящие от  } (f^{(\nu-1)},g^{(\nu-1)}), (f^{(\nu-2)},g^{(\nu-2)}),\dots = 0.$$
Т.е. эта система линейных соотношений обладает свойством треугольности. Это позволяет оценивать ранг полной (возмущенной) системы $\mathcal{L}(f,g;\Phi,\Psi)=0$ через ранг модельной $L(f,g;\Phi)=0$. В силу чего получаем оценку на размерность. Теорема доказана.

\vspace{3ex}

{\bf Замечание 12:}(a)  Отмеченное в доказательстве свойство треугольности системы (\ref{TC}) позволяет сделать более точное утверждение. Пусть $\nu \in \mathbb{Z}$  обозначим через $\mathcal{V}^{\nu}$ подпространство пространства векторных полей, т.ч. их весовое разложение не содержит компонент веса меньше, чем  ${\nu}$.  Тогда можно утверждать, что
 $$ \forall \; \nu \quad \dim ({\rm aut} \, M_{0} \cap \mathcal{V}^{\nu}) \leq  \dim \,({\rm aut} \, Q_0  \cap \mathcal{V}^{\nu}) .$$
(b)  Рассуждение, приведенное в доказательстве применимо не только к полной алгебре автоморфизмов, но и к любой ее подалгебре. В частности, в качестве таких подалгебр можно взять $gs$ и $st_{-}$.
И тогда получаем
$$\dim \,gs( M_{0}) \leq  \dim \,gs(Q_0),   \quad \dim \,st_{-}( M_{0}) \leq  \dim \,st_{-}(Q_0).$$
Поскольку условием голоморфной однородности ростка является соотношение $\dim \,gs( M_{0})=\dim M =2n+K$,
то из первого неравенства следует, что любая модельная поверхность голоморфно однородного ростка -- голоморфно однородна. Это утверждение ранее было доказано в \cite{MS1}.\\
(c) Размерность алгебры автоморфизмов ростка не зависит от введенной нами  $\mu$-градуировки, поэтому неравенство
из теоремы 11 можно заменить на
$$\dim \,{\rm aut} \, M_{0} \leq  \min \dim \,{\rm aut} \, Q_0(\mu) \mbox{ по всем }  \mu.$$

\vspace{3ex}

Как и в \cite{VB20} введем для $\mu$-типа $m=((m_1,k_1), \dots,(m_q,k_q))$ понятие старшего веса, а именно,
положим $\lambda=\lambda(m)=m_q$. Чтобы иметь возможность говорить о $\mu$-типе поверхности в точке $\xi$ этой поверхности, введем обозначение $m(\xi)$,   тогда старший вес в точке $\xi$ -- это $\lambda(\xi)=\lambda(m(\xi))$.  Далее, определим следующие подмножества $Q$:\\
$$
Q^m{'} =\{\xi \in Q: m(\xi)=m{'} \}, \quad  Q^{\lambda{'}}= \{ \xi \in Q: \lambda(\xi)=\lambda{'}  \}
$$
Подмножество вещественного аффинного пространства называем {\it полуалгебраическим}, если оно задано условиями вида $\{R^1(x)=0, \; R^2(x) \neq 0\}$, где $R^1$ и $R^2$ -- два конечных набора вещественных полиномов.

\vspace{3ex}

{\bf Теорема 13:} \\
(a)  Пусть $\xi \in  Q$, тогда старший вес в $\xi$  не превосходит старшего веса в нуле, т.е. $\lambda(\xi) \leq \lambda(0)$.\\
(b)  Множества $Q^m{'}, \;\;  Q^{\lambda{'}}$ -- полуалгебраичны для любых $m{'}$ и $\lambda{'}$.\\
(c)  Множество значений функций $m(\xi)$ и $\lambda(\xi)$ -- конечно.\\
{\it Доказательство:} Пункт (a) следует из нашего рассмотрения процедуры построения приведенной формы уравнений в точке модельной поверхности. Пункт (b) -- очевиден. Утверждение пункта (c), касающегося функции $\lambda(\xi)$ следует из (a). Часть утверждения, касающаяся $m(\xi)$,  следует из того, что для фиксированных старшем весе и коразмерности существует лишь конечное число $\mu$-типов $m'$.

\vspace{3ex}

{\bf Теорема 14:}  Пусть $\mu$-модельная поверхность $Q$ -- невырождена и голоморфно однородна,  тогда группа ее голоморфных автоморфизмов ${\rm Aut} \, Q$ -- это подгруппа группы бирациональных автоморфизмов $\mathbf{C}^N$ (группа Кремоны),  состоящая из отображений равномерно ограниченных степеней $d(\chi)$.  Константа, ограничивающая степени, зависит лишь от $N$
$$d(\chi) \leq C(N).$$
{\it Доказательство:}  Схема доказательства таких утверждений,  идущая от В.Каупа \cite{K},  многократно использовалась (см.\cite{VB21}).  По модулю того, что голоморфная однородность реализуется полиномиально-треугольными сдвигами, достаточно доказать это утверждение для элемента стабилизатора.  Для это требуется полиномиальность полей с оценкой степени  и наличие градуирующего поля.  Все это имеется в нашей ситуации. Таким образом, теорема доказана.

\vspace{3ex}

Вопрос о голоморфной однородности многообразия является весьма тонким. В работе \cite {VB21B} для модельных поверхностей (обычных, невзвешенных) был приведен простой критерий. Было показано, что точка модельной поверхности попадает в орбиту начала координат в том и только том случае, если ее Блум-Грэм-тип совпадает с Блум-Грэм-типом начала координат. Этот критерий остается в силе и для взвешенных модельных поверхностей.

\vspace{3ex}

{\bf Теорема 15:} Пусть $Q$ -- $\mu$-модельная поверхность и $\xi$ -- ее точка. Эта точка принадлежит орбите начала координат в группе автоморфизмов  $Q$ в том и только том случае, если $\mu$-тип $Q$ в точке $\xi$ по Блуму-Грэму-Степановой совпадает с  $\mu$-типом $Q$ в начале координат.\\
{\it Доказательство:}  Если $\mu$-тип $Q$ в точке $\xi$ совпадает с  $\mu$-типом $Q$ в начале координат, то автоморфизм
 строится так. Переносим начало координат пространства в точку $\xi$, переразлагаем определяющие полиномы по новым координатам и запускаем процесс построения приведенной формы новых уравнений. В результате уравнения принимают тот же самый вид (иначе изменится тип), что в начале координат. При этом получен треугольно-полиномиальный автоморфизм, переводящий начало координат в $\xi$. Обратно. Пусть имеется автоморфизм $\chi$, переводящий  начало координат в $\xi$. Его можно представить в виде композиции двух автоморфизмов $\chi =\chi_0 \circ \tilde{\chi}$, где $\chi_0$  -- это элемент стабилизатора нуля,  а $\tilde{\chi} \in GS$ -- элемент $G_{-}$, который тождественен на комплексной касательной. Тогда
 $\tilde{\chi}=\chi =(\chi_0)^{-1} \circ {\chi}$. Это автоморфизм, переводящий  начало координат в  $\xi$ с тождественным действием на комплексной касательной. Такие отображения сохраняют тип. Следовательно, $\mu$-типы в нуле и $\xi$  совпадают. Теорема доказана.

\vspace{3ex}

{\bf Следствие  16:}  $\mu$-модельная поверхность $Q$ голоморфно однородна тогда и только тогда, когда все ее точки имеют одинаковый $\mu$-тип (тип по Блуму-Грэму-Степановой с весом $\mu$).

\vspace{3ex}

{\bf 6. Примеры голоморфно однородных $\mu$-невырожденных модельных поверхностей}

\vspace{3ex}
В качестве {\bf первого примера} голоморфно однородной $\mu$-модельной поверхности можно предложить упомянутую во введении гиперповерхность \eqref{AL}.  В этом примере тип имеет следующий вид
\begin{eqnarray*}
\Bigl((\mu_1=1,  n_1=2),  (\mu_2=2,  n_2=1), \dots, (\mu_p=p,  n_p=1); \,(m_1=p+1,\, k_1=1) \Bigr),\\
 \mbox{где   }  p \geq 1, \; q=1, \; n=p+1, \; N=p+2. \qquad \qquad\qquad\qquad\qquad\qquad
\end{eqnarray*}

\vspace{3ex}

Для удобства читателя и в качестве иллюстрации общих рассмотрений, проведенных выше, приведем здесь описание автоморфизмов первой гиперповерхности из этой серии \cite{VB21} (теорема 15).  Пусть $p=2, \; n=3, \; N=4$.  Мы получаем гиперповерхность пространства $\mathbf{C}^4$
\begin{equation}\label{Q}
Q=\{v=2 \, {\rm Re} (z_1 \, \bar{\zeta} + z_2 \, \bar{\zeta}^2) \}
\end{equation}
Веса назначаем следующим образом
$$  [z_2]=[\zeta]=1, \; [z_1]=2, \;  [w]=[u]=3.$$
Для краткости векторное поле вида 
 $$X=2 \,{\rm Re} \left(  f _1\, \frac{\partial}{\partial z_1} + f _2\, \frac{\partial}{\partial z_2} +h\, \frac{\partial}{\partial \zeta} +g \, \frac{\partial}{\partial w} \right)$$
 будем записывать как $(f_1,f_2,h,g)$. 
Алгебра ${\rm Aut} \, Q$ имеет вид $g_{-3}+g_{-2}+g_{-1}+g_{0}+g_{1}$ причем
\begin{eqnarray} \label{DQ}
  g_{-3}=\{(\; 0, \quad 0, \quad 0, \quad d \;) \},   \qquad\qquad \qquad \qquad \qquad  \qquad \qquad \qquad \qquad \qquad \qquad \qquad \qquad \qquad \qquad \qquad \\
\nonumber
  g_{-2}=\{(\; a,\ \quad 0,\quad 0,\quad 2 \,i \, \bar{a}\, \zeta \; ) \},    \qquad \qquad\qquad \qquad \qquad \qquad \quad  \qquad \qquad \qquad \qquad \qquad \qquad \qquad \qquad \qquad\\
\nonumber
g_{-1}=\{(\; - 2 \, \bar{c} \, z_2 + i \, e \, \zeta , \quad  b,  \quad  c, \quad  2 \,i \, \bar{c}\, z_1+2 \, i \, \bar{b} \, \zeta^2 \;) \},   \qquad \qquad \qquad \qquad \quad \qquad \qquad \qquad \qquad \qquad \qquad\qquad\qquad\\
\nonumber
g_{0}= \{ ( \; \alpha_1 \,  z_1  - \bar{\alpha}_2 \,  \zeta^2 , \quad  (2 \, \alpha_2  - \alpha_3) \, z_2  + \alpha_2 \, \zeta, \quad
(\alpha_3-\bar{\alpha}_1)\, \zeta, \quad \alpha_3 \, w \; )  \},    \qquad \qquad \qquad\qquad \qquad \qquad \qquad \qquad \qquad \qquad \qquad\\
\nonumber
g_{1}=\{( \; 2\,i\,\bar{\beta}_1\,z_1\,\zeta+\beta_1\,w, \quad  2 \,i\, \bar{\beta}_1\, z_2\,\zeta -i\,\beta_1\, z_1+i\,\beta_2 \, \zeta^2, \quad i\,\bar{\beta}_1\,\zeta^2, \quad  2\,i\,\bar{\beta}_1\,\zeta w \; ) \},  \qquad \qquad \qquad  \qquad \qquad \quad  \qquad \qquad \qquad \qquad\\
a, \; b, \; c, \; \alpha_1, \; \alpha_2, \; \beta_1  \in \mathbf{C},   \quad d, \, e, \; \alpha_3, \; \beta_2 \in \mathbf{R}. \qquad \qquad \qquad \qquad\qquad \qquad \qquad \qquad\qquad \qquad \qquad \qquad
\end{eqnarray}
Запишем $g_{-1}$  в виде прямой суммы $g'_{-1}+st_{-}$, где
\begin{eqnarray*}
g'_{-1}=\{(- 2 \, \bar{c} \, z_2  ,\; b, \; c, \; 2 \,i \, \bar{c}\, z_1+2 \, i \, \bar{b} \, \zeta^2) \},\\
st_{-}=\{(  i \, e \, \zeta,\;0,\;0,\;0) \}.
\end{eqnarray*}
Тогда алгебра $gs$,  соответствующая группе "сдвигов", $\,$  $GS$ имеет вид $gs=g_{-3}+g_{-2}+g'_{-1}$.  Она параметризуется набором $(a,b,c,d)$, соответственно, $\dim \, gs=7$. Сама подгруппа $GS$, которая обеспечивает голоморфную однородность $Q$, состоит из преобразований вида
\begin{eqnarray}\label{QQ}
\nonumber
z_1 \rightarrow  A + z_1 ,  \qquad z_2 \rightarrow B + 2 \, \bar{A} \, \zeta + z_2, \qquad \zeta  \rightarrow  C + \zeta , \qquad\qquad \\
w    \rightarrow D  + 2 \, i \, (A \, \bar{B}  + A^2 \, \bar{C} +  (\bar{B} +2 \, A \, \bar{C})\, z_1 + \bar{A} \, z_2   +
\bar{A}^2 \, \zeta + \bar{C} \, z_1^2) + w,
\end{eqnarray}
где $(A,B,C,D)$ -- произвольная точка $Q$.

$\dim \, st_{-}=1$,  поле $(i \, \zeta,0,\;0,\;0)$  порождает группу $St_{-}$, которая имеет вид
$$z_1 \rightarrow z_1 +i \, t \, \zeta, \quad z_2 \rightarrow  z_2,  \quad  \zeta \rightarrow \zeta , \quad w  \rightarrow w.$$

Алгебра $g_0$ параметризуется набором $(\alpha_1,\alpha_2,\alpha_3)$, соответственно, $\dim \, g_0=5$.  Для вычисления группы $G_0$, соответствующей $g_0$, положим $\gamma=\alpha_1+\bar{\alpha}_1-\alpha_3$. Если $\gamma \neq 0$, то получаем
\begin{eqnarray*}
z_1  \rightarrow  \left( z_1 - \bar{\alpha}_2 \,  \left(\frac{e^{\gamma \, t}-1}{\gamma}\right) \, \zeta^2  \right) \, e^{\alpha_1 \,t}, \\
z_2  \rightarrow  \left( z_2 + \alpha_2 \, \left(\frac{e^{1-\bar{\gamma} \, t}}{\bar{\gamma}}\right) \, \zeta  \right) \, e^{(2 \,\alpha_1-\alpha_3)\,t},\\
\zeta \rightarrow \zeta \, e^{(\alpha_3-\bar{\alpha}_1)\,t}, \quad w \rightarrow w \, e^{\alpha_3 \,t}. \qquad\qquad
\end{eqnarray*}
Вырожденные направления $\gamma=0$ получаем предельным переходом.

Подалгебра $g_{+}$ состоит из единственной компоненты $g_1$, которая параметризуется набором $(\beta_1,\beta_2)$,  соответственно, $\dim \, g_1=3$.
Поле   $(0, \; i \, \zeta^2,\;0,\;0)$  из $g_{1}$ ($\beta_1=0, \; \beta_2=1$) порождает преобразование
\begin{equation}\label{G11}
z_1 \rightarrow  z_1, \;z_2 \rightarrow z_2 +i \, t \, \zeta^2, \;   \zeta \rightarrow \zeta , \; w  \rightarrow w.
\end{equation}
Преобразования из $g_1$ при $\beta_2=0$ имеют вид
\begin{eqnarray}\label{G12}
\nonumber z_1 \rightarrow\frac{z_1}{(1-i\, \bar{\beta}_1 \, \zeta \, t)^2},   \quad z_2   \rightarrow\frac{z_2-i\,\beta_1\, z_1\,t}{(1-i\, \bar{\beta}_1 \, \zeta \, t)^2},\\
\zeta \rightarrow \frac{\zeta}{1-i\, \bar{\beta}_1 \, \zeta \, t}, \quad    w \rightarrow \frac{w}{(1-i\, \bar{\beta}_1 \, \zeta \, t)^2}.
\end{eqnarray}
Преобразования (\ref{G11}) и (\ref{G12}) порождают группу $G_{+}$, соответствующую $g_{+}$.

\vspace{3ex}
В качестве {\bf второго примера} голоморфно однородной $\mu$-модельной поверхности можно предложить
гиперповерхность из работы \cite{ZS}.  Это гиперповерхность пространства  ${\bf C}^{n+1}$  с координатами $(z_1,...,z_n, \, w=u+i\,v)$, заданная уравнением
  $$v = 2 \; {\rm Re} (z_1 \bar{z}_2 + z_1^2  \bar{z}_3) +\sum_4^n  \pm |z_j|^2 .$$
  В этом примере тип имеет следующий вид
\begin{eqnarray*}
\Bigl((\mu_1=2,  n_1=2),  (\mu_2=4,  n_2=1), (\mu_3=3, n_3= n-3); \,(m_1=6, \, k_1=1) \Bigr),\\
\end{eqnarray*}
Эта гиперповерхность голоморфно однородна, 2-невырождена, размерность ее группы автоморфизмов равна $n^2+7$, причем это максимум в классе таких гиперповерхностей \cite{ZS}.
\vspace{3ex}

{\bf Третий пример}. Гиперповерхность из работы \cite{Kr}.  Это гиперповерхность пространства  ${\bf C}^{n+1}$  с координатами $(z_1,...,z_n, \, w=u+i\,v)$, заданная уравнением
\begin{equation}\label{4}
v = 2 \; {\rm Re} (z_1 \bar{z}_2)  +|z_1|^4 +  \sum_3^n  \pm |z_j|^2.
\end{equation}
    В этом примере тип имеет следующий вид
\begin{eqnarray*}
\Bigl((\mu_1=1,  n_1=1),  (\mu_2=3,  n_2=1), (\mu_3=2,  n_3=n-2); \,(m_1=4, \, k_1=1) \Bigr),\\
\end{eqnarray*}
Эта гиперповерхность голоморфно однородна и 1-невырождена (невырождена по Леви). В классе таких гиперповерхностей максимальную размерность $(n+2)^2-1$ имеет группа невырожденной стандартной гиперквадрики  $\{v=<z,\bar{z}>\}$, где  $<z,\bar{z}>$ -- невырожденная эрмитова форма. Гиперповерхность (\ref{4}) субмаксимальна, т.е. реализует следующее, после гиперквадрики, максимальное значение размерности группы, равное $n^2+4$.

\vspace{3ex}

{\bf Четвертый пример}. Семейство гиперповерхностей из списка работы \cite{Lb}.  Это гиперповерхности пространства  ${\bf C}^{3}$  с координатами $(z_1,z_2, \, w=u+i\,v)$, заданные уравнениями  $v = 2 \; {\rm Re} (z_1) \, {\rm Re} (z_2) + {\rm Re} (z_1) ^m$  ($m$ -- произвольное натуральное число).
\begin{eqnarray*}
\Bigl((\mu_1=1,  n_1=1),  (\mu_2=m-1,  n_2=1);  \,(m_1=m, \, k_1=1) \Bigr),\\
\end{eqnarray*}
Эти трубчатые  гиперповерхности голоморфно однородны и невырождены по Леви.  Гиперповерхности из примера три и четыре -- это некоторые обобщения гиперповерхности Винкельмана  $v=2 \,  {\rm Re} (z_1 \, \bar{z}_2) +|z_1|^4$.

\vspace{3ex}

Все рассмотренные здесь примеры -- это гиперповерхности. Получать на основе гиперповерхностей примеры большей коразмерности можно, рассматривая их прямое произведение. Но это далеко не единственная возможность.

\vspace{3ex}

Поскольку все эти $\mu$-модельные гиперповерхности невырождены и голоморфно однородны, то к ним применима теорема 11, т.е. их автоморфизмы бирациональны.

\vspace{3ex}

{\bf 7. Вполне $\mu$-невырожденные модельные поверхности}

$\ \ \ $  {\bf  и теория Танаки}

\vspace{3ex}

Также можно предложить достаточно широкий класс примеров однородных модельных гиперповерхностей для произвольного веса $\mu$. Эти $\mu$-модельные поверхности являются естественными обобщениями обычных вполне невырожденных модельных поверхностей \cite{VB04}.

Зафиксируем произвольную градуировку переменной $z$
$$\mu =((\mu_1, n_1), (\mu_2,  n_2), \dots, (\mu_p, \,n_p))$$
и рассмотрим последовательность пространств приведенных многочленов $\mathcal{N}_{\nu}$ всех натуральных весов. При этом не утверждается, что все $\kappa_{\nu}=\dim \, \mathcal{N}_{\nu} >0$. Перенумеруем ненулевые пространства $\mathcal{N}$. Будем нумеровать их не весом, а номером шага, на котором они появляются (т.е. пропускаем нулевые пространства).

Ясно что первый (минимальный) возможный вес для ненулевого пространства -- это $2 \,\mu_1$. При этом $\mathcal{N}_{1}$ -- это пространство эрмитовых форм от $z_1$. Размерность этого пространства $\kappa_{1}=n_1^2$.  Соответственно, вводим переменную $w_1$ и полагаем  $m_1=[w_1]=2 \,\mu_1$ и $k_1=\kappa_{т_1} =n_1^2$.

Далее среди приведенных полиномов от $z_1, \, z_2, \, u_1$ следует выбрать подпространство полиномов минимального,  после $2 \,\mu_1$ веса.  Если $\mu_2 < 2 \, \mu_1$, то это вес равный $\mu_1+\mu_2$, если же -- нет, то $3 \,\mu_1$.
После чего вводим переменную $w_2$ и полагаем  $m_2=[w_2]$ равными второму весу и $k_2=\kappa_{т_2}$. И так далее.
Останавливаемся на произвольном $q$-м шаге.  Теперь для $\beta=1,\dots,q-1$ положим
$\Phi_{\beta}$ -- равным набору базисных форм пространства $\mathcal{N}_{\beta}$, а $\Phi_{q}$ -- произвольному линейно независимому набору $(\Phi^1_{q},\dots,\Phi^k_{q})$ элементов $\mathcal{N}_{q}$, т.е. $1 \leq k \leq k_q$.
В результате в пространстве $\mathbf{C}^{n+K}, \; n=\sum_1^p \, n_{\alpha}, \; K=(\sum_1^{q-1} \,k_{\beta} +k)$  получаем поверхность
\begin{equation}\label{MTND}
Q=\{v_{\beta}=\Phi_{\beta}(z,\bar{z},u), \;\; \beta = 1, \dots, q\}
\end{equation}
Заметим, что произвол в выборе такой поверхности -- это произвол в выборе $q$ и набора координат формы старшего веса
$\Phi_q$.
\vspace{3ex}

{\bf Утверждение 17:}\\
(a) Любая такая поверхность $Q$ имеет конечный $\mu$-тип.\\
(b) Любая такая поверхность  голоморфно однородна.\\
(с) Если $m_q \geq 2 \, \mu_p+1$, то $Q$  голоморфно невырождена.\\
(d) Две такие поверхности, у которых совпадает параметр $q$ и старшие формы $\Phi_q$  линейно эквивалентны.\\
(e) Если кратность старшего веса $k_q$ принимает максимальное значение  $k_q = \dim \mathcal{N}_q$, то все такие
модельные поверхности линейно эквивалентны (т.е. такая поверхность единственна).\\
{\it Доказательство:}  (a) выполнен, т.к. $Q$ по построению имеет конечный $\mu$-тип в нуле. (b) -- т.к. переразложение уравнений в произвольной точке треугольно-полиномиальными преобразованиями приводится к исходному виду. (c) -- это следствие того, что при этом условии пространство $\mathcal{N}_{1}+\dots+\mathcal{N}_{q-1}$ содержит пространство всех
эрмитовых форм от $z$. (d) -- одна такая поверхность переводится в другую преобразованием
$$z  \rightarrow  z,  \quad  w_{\beta} \rightarrow \rho _{\beta} \, w_{\beta}, \;\; \beta = 1, \dots, q-1,   \;\;  w_{q} \rightarrow  w_{q},$$
где $\rho _{\beta}$ -- это вещественное невырожденное линейное преобразование, осуществляющее переход от одного базиса к другому в соответствующем пространстве. (e) Если $k_q = \dim \mathcal{N}_q$, то координаты $\Phi_q$ -- это базис  $\mathcal{N}_q$,  а два базиса связаны невырожденным линейным преобразованием. Утверждение доказано.

\vspace{3ex}

{\bf Определение 18:} Если построенная выше поверхность $Q$ является голоморфно невырожденной, то мы называем ее {\it вполне $\mu$-невырожденной} модельной поверхностью старшего веса $m_q$.

\vspace{3ex}
К любой  вполне $\mu$-невырожденной модельной поверхности применима теорема 12,  т.е. группа ее голоморфных автоморфизмов состоит из бирациональных преобразований равномерно ограниченных степеней.

\vspace{3ex}

Обычные, невзвешенные, вполне невырожденные модельные поверхности обладают следующим свойством.
Если старшая степень модельной вполне невырожденной поверхности $Q$  больше двух, то $g_{+}=0$.

Вполне невырожденные модельные поверхности, таким образом, делятся на два класса: невырожденные квадрики,
для которых  старшая степень $l=2$ и у которых  $g_{+}$ может быть нетривиальной, и остальные, у которых $l \geq 3$
и нет компоненты  $g_{+}$. Между этими двумя классами имеется еще одно важное различие. 
Критерием конечномерности алгебры квадрики являются два условия: линейная независимость эрмитовых форм и отсутствие у них общего ядра.  Это простые условия, в которые для квадрики превращаются два общих требования:  конечность Блум-Грэм-типа и голоморфная невырожденность.  При этом из одной конечности типа голоморфная невырожденность не следует. Тогда как при $l \geq3$
требование  вполне невырожденности -- это условие только на Блум-Грэм-тип. Голоморфная невырожденность  следует из этого условия на тип.

Для нашего взвешенного аналога полной невырожденности также можно ввести деление на такие классы.  Рассмотрим
последовательность $\{\mathcal{Q}_s\}, \; s=2,3,\dots$, где  $\mathcal{Q}_s$ -- последняя (т.е. такая, у которой все кратности максимальны) модельная вполне $\mu$-невырожденная поверхность старшего веса $s$.

\vspace{3ex}

{\bf Лемма 19:} Существует $\mathbf{s}$, т.ч. все  $\mathcal{Q}_s$ при $s \geq \mathbf{s}$ голоморфно невырождены.\\
{\it Доказательство:}  Начиная с некоторого $s$ все пространства $\mathcal{N}_s$ содержат пространство всех эрмитовых форм от $z$. Поэтому если кратность веса $s$ равна размерности $\mathcal{N}_s$, то среди  правых частей уравнений, определяющих $\mathcal{Q}_s$, содержится базис пространства эрмитовых форм. Это гарантирует голоморфную невырожденность. Лемма доказана.

\vspace{3ex}

{\bf Определение 20:} Минимальный из весов $\mathbf{s}$, существование которых доказано в лемме, назовем
{\it критическим} весом.

\vspace{1ex}
Ясно, что критический вес зависит от выбора базового набора весов $\mu$, т.е. $\mathbf{s}=\mathbf{s}(\mu)$.

\vspace{3ex}

Теперь мы готовы сформулировать взвешенный аналог $g_{+}$-гипотезы.

{\bf  Взвешенная $ g_{+}$-гипотеза:}  Пусть  $Q$  -- вполне $\mu$-невырожденная модельная поверхность старшего веса $l=m_q$, т.ч. $l > \mathbf{s}(\mu)$, тогда  $g_{+}=0$.

\vspace{3ex}

В доказательстве  стандартной  $g_{+}$-гипотезы  \cite{GK06}, \cite{SS18}, \cite{G18} значительную роль играет теория Н.Танаки. Центральные понятия этой теории -- это фундаментальная градуированная алгебра,  продолжение по Танаке и стандартная модель.  Мостом между теорией голоморфно однородных модельных поверхностей и теорией Танаки является алгебра  $g_{-}$, которая в теории Танаки является фундаментальной.

\vspace{3ex}

Непосредственно теория Танаки к взвешенным модельным поверхностям не применима.  В этой новой ситуации $g_{-}$ не является фундаментальной.  Однако, после некоторого редактирования теории Танаки, связь с теорией
голоморфно однородных взвешенных модельных поверхностей восстанавливается.

Как надо изменить теорию Танаки?  Вот небольшой эскиз этих изменений.

Первое изменение касается понятия фундаментальной градуированной алгебры. Пусть  $\mathfrak{g}=\sum_{-l}^{-1} \,\mathfrak{g}_j$ -- конечномерная градуированная алгебра Ли.  И пусть она порождается (как алгебра Ли) следующим набором своих компонент
$$\mathfrak{g}_{-\mu_1}, \; \mathfrak{g}_{-\mu_2}, \dots,  \mathfrak{g}_{-\mu_p}.$$
При этом не предполагается, что все $\mathfrak{g}_{-\nu} \neq 0$ для всех $1 \leq \nu \leq l$.
Такую алгебру мы будем называть {\it $\mu$-фундаментальной градуированной алгеброй Ли}.
Если при этом $\mathfrak{g}_{-l} \neq 0$,  то мы говорим, что это алгебра старшего веса $l$.  Ясно, что такая алгебра -- это аналог нашей подалгебры $gs$.

Мы будем говорить, что такая  $\mu$-фундаментальная алгебра  $\mathfrak{g}$ -- {\it невырождена}, если из того, что $Y \in g$ коммутирует со всеми образующими $\mathfrak{g}_{-\mu_1}, \; \mathfrak{g}_{-\mu_2}, \dots,  \mathfrak{g}_{-\mu_p}$,  следует, что $Y=0$.

Что касается аналога танаковского продолжения, то здесь появляется новый феномен. Это возможность продолжать $\mathfrak{g}$, присоединяя новые отрицательные компоненты. С алгебраической точки зрения каждое такое поле задает эндоморфизм алгебры $\mathfrak{g}$. Поэтому эти, отрицательные, весовые компоненты продолжения строятся как пространства эндоморфизмов с необходимым набором условий (тождество Якоби).  Так мы получаем набор новых отрицательных компонент $\mathfrak{st}_{-}$. Это аналог нашей компоненты $st_{-}$. После этого, так же, как и раньше, мы можем определить компоненту $\mathfrak{g}_0$ как пространство дифференцирований на $\mathfrak{g}_{-}=\mathfrak{g}+\mathfrak{st}_{-}$  (аналог нашего $g_{-}$). Далее, как и в невзвешенном случае, рекуррентно строится  последовательность компонент положительных весов, как пространств операторов на уже построенных компоненентах с условием, имитирующим тождество Якоби. Т.е. мы получаем $\mathfrak{g}_{+}=\mathfrak{g}_1 + \mathfrak{g}_2 + \dots$  (аналог нашего $g_{+}$).

  Можно надеяться,  однако это нуждается в доказательстве, что если $\mathfrak{g}$ -- $\mu$-фундаментальна и невырождена, то:\\
(a) существует $ j>0$ такое, что $\mathfrak{g}_{j}=0$.\\
(b) если для некоторого $ j>0$ оказалось, что $\mathfrak{g}_{j}=0$, то  $\mathfrak{g}_{j+1}=0$.
Т.е. в таком случае  продолжение алгебры $\mathfrak{g}$ конечно градуировано и конечномерно. \\
(c) аналог стандартной модели $\mathfrak{G_(\mathfrak{g})}$ как $CR$-многообразие эквивалентен $G_{-}$ для голоморфно однородной взвешенной модельной поверхности, у которой $gs=\mathfrak{g}$.

Если мы планируем использовать эту технику для наших целей, то нам потребуется снабдить образующие $\mathfrak{g}_{-\mu_1}, \; \mathfrak{g}_{-\mu_2}, \dots,  \mathfrak{g}_{-\mu_p}$ комплексной структурой $J$. В этом случае мы будем говорить, что мы имеем $CR$-фундаментальную $\mu$-градуированную алгебру Ли.

Аналогичным образом мы можем построить стандартную модель -- голоморфно однородное $CR$-многообразие, Блум-Грэм-тип которого кодируется алгеброй $\mathfrak{g}$. В частности,  $\mathfrak{g}_{-\mu_1}+ \mathfrak{g}_{-\mu_2}+ \dots +  \mathfrak{g}_{-\mu_p}$ -- это его комплексная касательная. Для этого следует рассмотреть связную группу Ли $\mathfrak{G}$, порожденную $\mathfrak{g}$ как вещественное подмногообразие в комплексной группе Ли $\mathfrak{G}^c$, порожденной $\mathfrak{g}^c$ -- комплексификацией $\mathfrak{g}$.

Вполне аналогично  понятию $CR$-универсальной фундаментальной градуированной алгебры Ли  вводим понятие   $CR$-универсальной $\mu$-фундаментальной градуированной алгебры Ли. Это такая фундаментальная градуированная алгебра Ли $\mathfrak{g}$, что в процессе ее порождения  базовым набором $\{\mathfrak{g}_{-\mu_j}\}$ при каждом коммутировании двух компонент рост размерности максимален.

Вопрос о тривиальности компонент положительного веса для универсальной взвешенной алгебры -- это алгебраическая версия сформулированной выше взвешенной $g_{+}$-гипотезы для взвешенных вполне невырожденных модельных поверхностей.

{\bf Взвешенная $\mathfrak{g}_{+}$-гипотеза } (алгебраическая версия):   Пусть $\mathfrak{g}$ --  невырожденная $CR$-универсальная $\mu$-фундаментальная алгебра Ли старшего веса $l$, т.ч. $l >\mathbf{s}$, тогда  $\mathfrak{g}_{+}=0$ (нет продолжений положительного веса).

Определение критического веса апеллирует к внешним по отношению к теории градуированных алгебр понятиям. Чтобы избежать этого можно последнее условие заменить на достаточное неравенство $l > 2 \mu_q$.

\vspace{3ex}

{\bf 8. Мультивесовая техника и оценка размерности алгебры автоморфизмов}

\vspace{3ex}

Модельные поверхности интересны, прежде всего тем, что каждая модельная поверхность является самой голоморфно симметричной по отношению к своим возмущениям (см. теорему 5).  Мультивесовая техника, т.е. использование различных весов в задаче простроения оценки размерности группы локальных голоморфных автоморфизмов, открывает весьма широкие перспективы.

\vspace{3ex}

Пусть  $M$ -- вещественно аналитическое голоморфно невырожденное порождающее $CR$-подмногообразие комплексного линейного пространства $\mathbf{C}^N$. Если $M$ имеет всюду бесконечный тип по Блуму-Грэму, то вопрос о размерности
локальных групп автоморфизмов ростков такого многообразия был недавно и весьма детально разобран в работе М.Степановой \cite{MS2}.  В точке общего положения такая локальная группа имеет размерность либо ноль, либо бесконечность.  На особом собственном аналитическом подмножестве возможна и положительная размерность.

Поэтому мы можем ограничиться рассмотрением случая многообразия, имеющего в точке общего положения конечный тип.
Причем в таком случае для задачи оценивания размерности мы можем ограничиться оценкой в точке конечного типа. Эта оценка, очевидно, будет выполнена и на особом подмножестве, которое является собственным аналитическим подмножеством.

Итак, пусть  $M_{\xi}$ -- росток многообразия $M$ и пусть $n$ -- его $CR$-размерность,  $K$ -- коразмерность. Тогда мы можем выбрать локальные координаты в точке $\xi$ так, что уравнение ростка примет вид ${\rm Im} \, w = F(z,\bar{z},{\rm Re}  \, w)$, где $(z \in \mathbf{C}^n,\; w=u+i\,v \in \mathbf{C}^K)$ -- координаты объемлющего пространства и вещественно аналитическая вектор-функция $F$ и ее первые производные обращаются в ноль в начале координат.

Разобьем произвольно число $n$ в сумму положительных слагаемых $n=n_1+\dots+n_p$. Соответственно, возникнет
разложение $\mathbf{C}^n$ в прямую сумму слагаемых вида $\mathbf{C }^{n_{\alpha}}$.  Назначим, тоже произвольно,
веса для элементов каждого прямого слагаемого, пусть $[z_{\alpha}]=\mu_{\alpha}$, где $ (\mu_{1}, \dots, \mu_{p})$ -- возрастающая последовательность натуральных чисел. Это формирует $\mu=((\mu_1,n_1), \dots,(\mu_p,n_p))$ -- весовое разложение переменной $z$. В соответствии с процедурой, описанной выше, среди переменных группы $w$ тоже возникнет весовое разложение $w=(w_1, \dots, w_q)$ и, в соответствии с типом по Блуму-Грэму-Степановой, будут назначены веса и кратности.  После некоторой полиномиальной замены координат уравнения ростка можно будет записать в приведенном виде и получить модельную поверхность $Q$. Если $Q'$ еще одна модельная поверхность того же ростка при том же выборе разбиения на $n_{\alpha}$ (кратности) и весов $\mu_{\alpha}$, то $Q'$ эквивалентна $Q$. Отображение осуществляется
обратимым квазилинейным отображением (утверждение 2, (b)). В этом смысле такая модельная поверхность единственна. При фиксированном ростке $M_{\xi}$ она зависит только от $\mu$.  Обозначим ее через $Q(\mu)$.

Препятствием для получения оценки является  голоморфная вырожденность модельной поверхности при любом выборе веса $\mu$.   Рассмотрим при фиксированном ростке конечного типа $M_{\xi}$ счетную совокупность взвешенных модельных
поверхностей $\mathcal{Q}=\{Q(\mu,M_{\xi})\}$ для всех возможных $\mu$.

\vspace{3ex}

{\bf  Определение 21:}    Росток  $M_{\xi}$ назовем {\it правильным}, если существует такой вес $\mu$, что соответствующая  модельная поверхность $Q(\mu,M_{\xi})$ --  невырождена (конечный тип плюс голоморфная невырожденность). Совокупность таких весов  $\mathfrak{M}(M_{\xi})$, которое в этом случае не пусто, назовем множеством {\it правильных весов}.

\vspace{3ex}
Ясно что правильный росток имеет конечный тип. Для правильных ростков из теоремы 10 сразу получаем .

\vspace{3ex}

{\bf  Утверждение 22:} Пусть $M_{\xi}$ -- правильный росток, тогда
$$ \dim \, {\rm aut} \,  M_{\xi} \leq  \min  \dim \, {\rm aut} \, Q_0(\mu,M_{\xi})   \mbox{ по всем  правильным весам}.$$

\vspace{3ex}
Заметим, хотя это и не снимает всех вопросов, что росток общего положения -- правильный.

\vspace{3ex}

В качестве примера неправильной гиперповерхности в $\mathbf{C}^3$ можно предложить хорошо известный  $"$световой конус$"$, который в координатах $z_j=x_j+i \,y_j, \; j=1,2,3$, задается уравнением $y_3^2=y_1^2+y_2^2, \; y_3>0$.

Описываемый нами подход непосредственно не применим к неправильным многообразиям и их росткам. Они требуют специального подхода.  В качестве такого подхода в работах \cite{VB05} и \cite{VB21} была предложена процедура оценивания, основанная на модифицированной конструкции Пуанкаре (рекурсия на глубину больше, чем один). Однако эта процедура является технически более сложной.

\vspace{3ex}

Использование модельных поверхностей при изучении автоморфизмов ростка имеет очевидную мотивировку. Росток  $M_{\xi}$ -- аналитический объект, а его модельная поверхность $Q(M_{\xi},\mu^0)$ -- алгебраический.  Модельная поверхность и ее автоморфизмы проще исходного ростка и его автоморфизмов. При этом все же следует отметить, что если размерность пространства велика, то полиномиальные квазиоднородные формы $\Phi$,  которые задают уравнения модельной поверхности, -- это полиномы большого числа переменных и общая модельная поверхность недоступна для непосредственного анализа. Однако, меняя начальный правильный вес $\mu^0$ на какой-либо другой правильный вес $\mu^1$ (для $Q^0$), мы можем сопоставить старой модельной  поверхности $Q^0=Q(M_{\xi},\mu^0)$ новую модельную поверхность $Q^1=Q(Q^0,\mu^1)$. Применив еще раз теорему 18,  мы получим оценку размерности автоморфизмов исходного ростка $M_{\xi}$ через размерность автоморфизмов $Q^1$.  Эту операцию можно повторять, уменьшая число ненулевых мономов в текущих квазиоднородных формах $\Phi$.  При этом модельная поверхность становится проще, однако полученная оценка может ухудшаться. Процесс закончится, как только мы получим модельную поверхность, для которой  нет следующего правильного веса. Отметим, что на каждом шаге, начиная с нулевого, мы вправе выбирать произвольный правильный вес по отношению стартовому ростку или же к текущей модельной поверхности.  Таким образом, таких цепочек может быть очень много. Каждая такая цепочка за конечное число шагов заканчивается своим терминальным звеном $Q^{\infty}$, полученым с помощью выбора веса ${\mu}^{\infty}$. И в качестве оценки автоморфизмов исходного ростка мы можем предложить размерность автоморфизмов $Q^{\infty}$. Эта последняя в цепочке модельная поверхность является функцией стартового ростка $M_{\xi}$ и цепочки весов
$$ \mathfrak{Z}=(\mu^0 \rightarrow \mu^1 \rightarrow \dots \rightarrow {\mu}^{\infty}).$$
В этом случае мы говорим, что цепочка $\mathfrak{Z}$ является правильной цепочкой для $M_{\xi}$ и что $Q^{\infty}=Q^{\infty}(M_{\xi},\mathfrak{Z})$. При этом мы можем написать
$$ \dim \, {\rm aut} \,  M_{\xi} \leq    \dim \, {\rm aut} \, Q_0^{\infty}(M_{\xi},\mathfrak{Z}) $$

\vspace{3ex}

Рассмотрим внимательнее случай гиперповерхности в пространстве $\mathbf{C}^{n+1}$ с координатами $z=(z_1,\dots,z_n), \; w=u+i\,v$.  Росток гиперповерхности $\Gamma_{0}$ -- это график вещественно аналитической функции вида
$$\{v=F(z,\bar{z}) +u \,G(z,\bar{z},u) \},  \; F(z,\bar{z})=\sum \, c_{\alpha \beta} \, z^{\alpha}\, \bar{z}^{\beta}, \;\alpha \in \mathbf{Z}_{+}^n, \; \beta  \in \mathbf{Z}_{+}^n.$$
Коэффициенты $c_{\alpha \beta}$ удовлетворяют условиям сходимости ряда и его вещественности.  Мы также можем предполагать, что разложение $F$ не содержит плюригармонических слагаемых, т.е. $F(z,0)=0$.  Функции $F$ можно также сопоставить  носитель ряда  -- совокупность узлов решетки $\mathbf{Z}_{+}^{2 \,n}$, для которых соответствующий моном стоит с ненулевым коэффициентом и его выпуклую оболочку $\mathfrak{N}$ -- многогранник Ньютона.

   На нулевом шаге мы выбираем правильный вес $\mu^0=(\mu^0_1,\dots,\mu^0_n)$, т.е. полагаем, что веса $z_{\nu}$ и $\bar{z}_{\nu}$ равны $\mu^0_{\nu}$, и получаем невырожденную полиномиальную модельную поверхность
   $$Q^0=\{v=P^0(z,\bar{z})=\sum \, p^0_{\alpha \beta} \, z^{\alpha} \, \bar{z}^{\beta} \},$$
    где $P^0$ - это квазиоднородный полином от  $(z ,\bar{z})$ некоторого веса $m^0 \geq 2$ без плюригармонических членов.  Таким образом мультииндексы
   $(\alpha,\beta) \in \mathbf{Z}_{+}^{2n}$ связаны линейным соотношением
$$
\mu^0_1 \,(\alpha_1+\beta_1)+\dots+\mu^0_n \,(\alpha_n+\beta_n)=m^0.
$$
Это соотношение выделяет в пространстве мономов некоторое конечномерное подпространство, пусть $\mathfrak{N}^0$ -- это многогранник Ньютона  $P^0$.  Следующий вес $\mu^1=(\mu^1_1,\dots,\mu^1_n)$ дает нам нам новую градуировку, которая с помощью соотношения
$$
\mu^1_1 \,(\alpha_1+\beta_1)+\dots+\mu^1_n \,(\alpha_n+\beta_n)=m^1.
$$
высекает на многограннике  $\mathfrak{N}^0$ некоторую грань $\mathfrak{N}^1$ меньшей размерности. И так далее до терминального веса $\mu^{\infty}$.  При этом мы получаем модельную поверхность $Q^{\infty}$ и ее многогранник $\mathfrak{N}^{\infty}$, которые для этой цепочки дают максимальное упрощение. Ясно, что число шагов ограничено размерностью пространства, т.е. числом $2n$.

\vspace{3ex}

Наш алгоритм работает в направлении уменьшения сложности модельной поверхности и останавливается при угрозе получения голоморфно вырожденной модельной поверхности.  Можно предложить альтернативный подход, который начинает с наипростейшей голоморфно вырожденной поверхности и, двигаясь в сторону усложнения, останавливается по достижении голоморфной невырожденности.  Начинаем с произвольной вершины $V_0$ на границе исходного многогранника $\mathfrak{N}$ (строго говоря -- это пара симметричных вершин). Стартовая модельная поверхность -- это график вещественной части этого монома.  Вес $\mu^0$ выбираем так, что опорная гиперплоскость, проходящая через $V_0$, не имеет других пересечений с $\mathfrak{N}$.  Проверку на голоморфную невырожденность осуществляем как проверку на конечную невырожденность. Если поверхность голоморфно вырождена, то выбираем второй моном -- вершину $V_1$ на границе $\mathfrak{N}$, смежную с $V_0$.  Вес $\mu^1$ выбираем так, что опорная гиперплоскость, проходящая через ребро $[V_0,V_1]$, не имеет других пересечений с $\mathfrak{N}$. Проверка на голоморфную невырожденность. И так далее до достижения голоморфной невырожденности.  Следует отметить, что модельная поверхность, чье уравнение связывает не все координаты, не может быть голоморфно невырожденной.

\vspace{3ex}

Первый описанный здесь алгоритм естественно назвать нисходящим, в второй -- восходящим.

\vspace{3ex}

Мономиальная модельная гиперповерхность -- это гиперповерхность вида $\{v = 2 \, {\rm Re} (z^{\alpha} \, \bar{z}^{\beta})\}$.
Если $z$ -- переменная  размерности $n=1$ или $n=2$, то такая гиперповерхность может быть голоморфно невырожденной. Вот примеры.
\begin{eqnarray*}
\{v=|z|^2 \} \quad \mbox{ для } n=1,\\
\{v = 2 \, {\rm Re} (z_1 \, \bar{z}_2)\} \quad \mbox{ для } n=2.
\end{eqnarray*}
Если же $n \geq 3$ -- это невозможно.

\vspace{3ex}

{\bf  Утверждение 23:} Если $n \geq 3$, то гиперповерхность  $\Gamma=\{v = 2 \, {\rm Re} (z^{\alpha} \, \bar{z}^{\beta})\}$ голоморфно вырождена.\\
{\it Доказательство:}  Голоморфная вырожденность вещественно аналитического многообразия эквивалентна его голоморфной вырожденности в произвольной точке. Если $\alpha = \beta$, рассмотрим замену координат вида
$z_1 \rightarrow  z^{\alpha}$, если $\alpha \neq \beta$, замену $z_1 \rightarrow  z^{\alpha},  \; z_2 \rightarrow  z^{\beta}$ (остальные координаты -- на месте).  В точке общего положения эта голоморфная замена локально обратима. После этой замены уравнения гиперповерхности принимают вид $v=|z_1|^2 $ или  $v = 2 \, {\rm Re} (z_1 \, \bar{z}_2)$. Уравнения содержат не все координаты, поэтому обе гиперповерхности и их прообраз $\Gamma$ голоморфно вырождены.\\

\vspace{1ex}

Применяя аналогичные соображения можно показать, что минимальное число мономов, обеспечивающее голоморфную невырожденность, равно $\{\frac{n+1}{2}\}$ ($\{x\}$ -- целая часть $x$).  В связи с этим можно предложить следующее усовершенствование второго предложенного (восходящего) алгоритма. Его работу надо начинать не с вершины многогранника, а с грани размерности $\{\frac{n+1}{2}\}$.

В работе этого алгоритма есть много степеней  свободы. Это свобода в выборе граней растущих размерностей или, что то же самое,  в выборе определяющих их весов.  Условием окончания работы алгоритма является достижение голоморфной невырожденности. Голоморфная невырожденность равносильна конечной невырожденности в точке общего положения.
Поэтому выбор грани и веса следует подчинить условию роста размерности совокупности производных градиента определяющего  полинома (см. определение конечной невырожденности).

Результатом работы этих алгоритмов, как нисходящего, так и восходящего, является {\it минимальная} модельная поверхность.
Т.е. такая голоморфно невырожденная поверхность, модельная по отношению к исходному ростку, т.ч. удаление любого ее монома делает ее голоморфно вырожденой.

\vspace{5ex}

{\bf 9. Вопросы на будущее}

\vspace{3ex}

В связи с описанным выше подходом к оценке размерности представляют интерес следующие вопросы.

\vspace{3ex}

{\bf  Вопрос 24:} Найти конструктивный критерий правильности ростка.

     Далее, возможно, имеются явные способы выделять из исходного многогранника  $\mathfrak{N}$ его модельное подмножество, а не получать его в качестве результата работы алгоритма.  Такой  способ представлял бы несомненный интерес.

\vspace{3ex}

Нельзя одновременно упрощать модельное подмножество и минимизировать получающуюся оценку размерности алгебры автоморфизмов.  Если мы отказываемся от простоты и ищем модельную поверхность с наименьшей размерностью, то,
как нетрудно понять, мы должны выбирать модельные подмножества  $\mathfrak{N}$, лежащие в граничных гипергранях.
Однако гиперграней несколько.

\vspace{1ex}

{\bf  Вопрос 25:} Как найти граничную гипергрань $\mathfrak{N}$, т.ч. соответствующая модельная поверхность дает минимальную размерность алгебры автоморфизмов?  Или, что то же самое, как найти соответствующий вес?

\vspace{3ex}

{\bf  Вопрос 26:} Взвешенная $g_{+}$-гипотеза, сформулированная в пункте 5.  Гипотеза имеет две редакции: геометрическую и алгебраическую.

\vspace{3ex}

И вообще, помещая теорию модельных поверхностей в новый взвешенный контекст, мы можем все старые вопросы рассмотреть  с весовой точки зрения.  Это имеет смысл по отношению как к доказанным утверждениям, так и к гипотезам (см. список в конце \cite{VB20}).

\end{document}